\documentclass[preprint,12pt]{elsarticle}
\usepackage{amsmath}
\usepackage{amssymb}

\input{tcilatex}

\begin{document}

\begin{frontmatter}



\title{
Equilibrium existence results for a class of discontinuous games
}

\author{Monica Patriche}

\address{
University of Bucharest
Faculty of Mathematics and Computer Science
    
14 Academiei Street
   
 010014 Bucharest, 
Romania
    
monica.patriche@yahoo.com }

\begin{abstract}
We introduce the notions of w-lower semicontinuous and almost w-lower semicontinuous 
correspondence with respect to a given set and prove a new fixed-point theorem. We also 
introduce the notion of correspondence with e-LSCS-property. As applications we obtain 
some new equilibrium theorems for abstract economies and for generalized multiobjective games.

\end{abstract}


\end{frontmatter}



\label{}





\bibliographystyle{elsarticle-num}
\bibliography{<your-bib-database>}







\section{Introduction}

In [16] W. Shafer and H. Sonnenschein proved the existence of equilibrium of
an economy with finite dimensional commodity space and irreflexive
preferences represented as correspondences with open graph. They generalized
the work of J. Nash [14], who first proved a theorem of equilibrium
existence for games where the player's preferences are representable by
continuous quasi-concave utilities and the work of G. Debreu, who proved the
existence of equilibrium in a generalized N-person game or an abstract
economy [3]. N. C. Yannelis and N. D. Prahbakar [20] developed new
techniques based on selection theorems and fixed-point theorems. Their main
result concerns the existence of equilibrium when the constraint and
preference correspondences have open lower sections. They worked within
different framework (countable infinite number of agents, infinite
dimensional strategy spaces). K. J. Arrow and G. Debreu proved the existence
of Walrasian equilibrium in [1]. To sum up, the significance of equilibrium
theory stems from the fact that it develops important tools to prove the
existence of equilibrium for different types of games.

A. Borglin and H. Keiding [2] used new concepts of K. F.-correspondences and
K. F.-majorized correspondences for their existence results. The second
notion was extended by Yannelis and Prabhakar [20] to L-majorized
correspondences. In [21], G. X. Yuan proposed a model of abstract economy
more general than that introduced by A. Borglin and H. Keiding in [2],
meaning that each constraint mapping has been divided into two parts, A and
B, because the set of the fixed points of \ the "small" correspondence may
not be rich enough.

Most of the existence theorems of equilibrium deal with preference
correspondences which have lower open sections or are majorized by
correspondences with lower open sections. Within the last years, some
existence results were obtained for lower semicontinuous and upper
semicontinuous correspondences. Some results concerning fixed point theorems
for lower semicontinuous correspondences or equilibrium existence for
economies with lower semicontinuous and Q-majorized correspondences can be
found in [4], [5], [10], [18]. E. Michael gave some selection theorems for
lower semicontinuous correspondences. His main results can be found in
[11]-[13].

In this paper, we define several types of correspondences: w-lower
semicontinuous and almost w-lower semicontinuous with respect to a given set
and also correspondences having e-LSCS-property. We prove a fixed point
theorem for almost w\textit{-}lower semicontinuous correspondences. This
result is a Wu like fixed point theorem [18]. We use this theorem to prove
the equilibrium existence for abstract economies which have w\textit{-}lower
semicontinuous\textit{\ }constraint correspondences. We use a technique of
approximation to prove an equilibrium existence theorem for correspondences
with e-LSCS-property.

We give slight generalizations of the equilibrium notions defined by W. K.
Kim and X. P. Ding in [9] and we also prove the existence of generalized
weighted Nash equilibrium and of generalized Pareto equilibrium for a
generalized multiobjective game having w\textit{-}lower semicontinuous
constraints.

The paper is organized in the following way: Section 2 contains
preliminaries and notation. The fixed point theorem is presented in Section
3 and the equilibrium theorems are stated in Section 4. Section 5 contains
the model of a constrained multiobjective game and a Pareto equilibrium
existence result.\bigskip

\section{\textbf{Preliminaries and notation\protect\smallskip }}

We shall denote by $\mathbb{R}_{+}^{m}:=\{u=(u_{1},u_{2},...,u_{m})\in 
\mathbb{R}^{m}:u_{j}\geq 0$ $\forall j=1,2,...,m\}$ and int$\mathbb{R}%
_{+}^{m}:=\{u=(u_{1},u_{2},...,u_{m})\in \mathbb{R}^{m}:u_{j}>0$ $\forall
j=1,2,...,m\}$\textit{\ }the non-negative othant of $\mathbb{R}^{m}$ and
respective the non-empty interior of $\mathbb{R}_{+}^{m}$ with the topology
induced in terms of convergence of vector with respect to the Euclidian
metric. For each $u,v\in \mathbb{R}^{m}$, $u\cdot v$ denote the standard
Euclidian inner product.\medskip

Now we present some notations and results concerning the theory of
correspondences.\medskip

Let $A$ be a subset of a topological space $X$. $\tciFourier (A)$ denotes
the family of all nonempty finite subset of A. $2^{A}$ denotes the family of
all subsets of $A$. cl$A$ denotes the closure of $A$ in $X$. If $A$ is a
subset of a vector space, co$A$ denotes the convex hull of $A$. If $F$, $G:$ 
$X\rightarrow 2^{Y}$ are correspondences, then co$G$, cl$G$, $G\cap F$ $:$ $%
X\rightarrow 2^{Y}$ are correspondences defined by $($co$G)(x)=$co$G(x)$, $($%
cl$G)(x)=$cl$G(x)$ and $(G\cap F)(x)=G(x)\cap F(x)$ for each $x\in X$,
respectively. The graph of $T:X\rightarrow 2^{Y}$ is the set Gr$%
(T)=\{(x,y)\in X\times Y\mid y\in T(x)\}.$

The correspondence $\overline{T}$ is defined by $\overline{T}(x)=\{y\in
Y:(x,y)\in $cl$_{X\times Y}$Gr$T\}$ (the set cl$_{X\times Y}$Gr$(T)$ is
called the adherence of the graph of $T$)$.$ It is easy to see that cl$%
T(x)\subset \overline{T}(x)$ for each $x\in X.$

\begin{remark}
$\overline{T}(x)=$cl$T(x)$ for each $x\in X$ if $T$ has a closed graph in $%
X\times Y$ (by Theorem 7.1.15 in [8], it follows that in particular, $T$ has
a closed graph when $Y$ is regular and cl$T$ is upper semicontinuous with
closed values).
\end{remark}

Let $X$, $Y$ be topological spaces and $T:X\rightarrow 2^{Y}$ be a
correspondence. $T$ is said to be \textit{lower semicontinuous} \textit{%
(l.s.c)} if for each x$\in X$ and each open set $V$ in $Y$ with $T(x)\cap
V\neq \emptyset $, there exists an open neighbourhood $U$ of $x$ in $X$ such
that $T(y)\cap V\neq \emptyset $ for each $y\in U$. $T:X\rightarrow 2^{Y}$
is said to be \textit{almost lower semicontinuous} if for each $x\in X$ and
each open set $V$ in $Y$ with $T(x)\cap V\neq \phi ,$ there exists an open
neighborhood $U$ of $x$ in $X$ such that $T(x)\cap \overline{V}\neq \phi $
for each $z\in U.$\medskip

For some known results about lower semicontinuity, which will be used in our
proofs, we refer the reader to [8].

\begin{proposition}
(Lemma 1 in [17], Proposition 2.5 in [11]$).$ Let $X$ and $Y$ be topological
spaces and let $T_{1},T_{2}$ be two l.s.c correspondences from $X$ to $Y.$
If $T_{1}$ has open values and $T_{1}(x)\cap T_{2}(x)\neq \emptyset ,$ for
every $x\in X,$ then, the correspondence $T$ defined by $T(x)=T_{1}(x)\cap
T_{2}(x)$ is l.s.c too.
\end{proposition}

\begin{proposition}
(Theorem 1.6, pag 25 in [21]). Let $X$ be a topological space, $E$ be a
topological vector space and $Y$ be a non-empty subset of $E.$ Suppose $%
S:X\rightarrow 2^{Y}$ is a lower semicontinuous correspondence and $V$ is
any nonempty open subset of $E.$ Then the correspondence $T:X\rightarrow
2^{Y}$ defined by $T(x)=(S(x)+V)\cap Y$ for each $x\in X$ has an open graph
in $X\times Y.$
\end{proposition}

We also need a version of Lemma 1.1 in [21]. For the reader's convenience,
we include its proof below.

\begin{lemma}
Let $X$ be a topological space, $Y$ be a nonempty subset of a locally convex
topological vector space $E$ and $T:X\rightarrow 2^{Y}$ be a correspondence$%
. $ Let \ss\ be a basis of neighbourhoods of $0$ in $E$ consisting of open
absolutely convex symmetric sets. Let $D$ be a compact subset of $Y$. If for
each $V\in $\ss , the correspondence $T^{V}:X\rightarrow 2^{Y}$ is defined
by $T^{V}(x)=(T(x)+V)\cap D$ for each $x\in X,$ then $\cap _{V\in \text{\ss }%
}\overline{T^{V}}(x)\subseteq \overline{T}(x)$ for every $x\in X.$
\end{lemma}

\textit{Proof.} Let be $x$ and $y$ be such that $y\in \cap _{V\in \text{\ss }%
}\overline{T^{V}}(x)$ and suppose, by way of contradiction, that $y\notin 
\overline{T}(x).$ This means that $(x,y)\notin $clGr$T,$ so that there
exists an open neighborhood $U$ of $x$ and $V\in $\ss\ such that:

$(U\times (y+V))\cap $Gr$T=\emptyset .\ \ \ \ \ \ \ \ \ \ \ \ \ \ \ \ \ \ \
\ \ \ \ \ \ \ \ \ \ \ \ \ (1)$

Choose $W\in $\ss\ such that $W-W\subseteq V$ (e.g. $W=\frac{1}{2}V)$. Since 
$y\in \overline{T^{W}}(x)$, then $(x,y)\in $clGr$T^{W},$ so that

\begin{equation*}
(U\times (y+W))\cap \text{Gr}T^{W}\neq \emptyset .\ \ \ \ \ \ \ \ \ \ \ \ \
\ \ \ \ \ \ \ \ \ \ \ \ \ \ \ \ \ \ \ 
\end{equation*}

There are some $x^{\prime }\in U$ and $w^{\prime }\in W$ such that $%
(x^{\prime },y+w^{\prime })\in $Gr$T^{W},$ i.e. $y+w^{\prime }\in
T^{W}(x^{\prime }).$ Then, $y+w^{\prime }\in D$ and $y+w^{\prime }=y^{\prime
}+w^{^{\prime \prime }}$ for some $y^{\prime }\in T(x^{\prime })$ and $%
w^{^{\prime \prime }}\in W.$ Hence, $y^{\prime }=y+(w^{\prime }-w^{^{\prime
\prime }})\in y+(W-W)\subseteq y+V,$ so that $T(x^{\prime })\cap (y+V)\neq
\emptyset .$ Since $x^{\prime }\in U,$ this means that $(U\times (y+V))\cap $%
Gr$T\neq \emptyset ,$ contradicting (1).\medskip

We present first Wu's Theorem 1 in [18], which will be generalized in the
next section.\medskip

\begin{theorem}[18]
\textit{Let }$I$\textit{\ be an index set. For each }$i\in I,$\textit{\ let }%
$X_{i}$\textit{\ be a nonempty convex subset of a Hausdorff locally convex
topological vector space }$E_{i}$\textit{, }$D_{i}$\textit{\ a non-empty
compact metrizable subset of }$X_{i}$\textit{\ and }$S_{i},T_{i}:X%
\rightarrow 2^{D_{i}}$\textit{\ two correspondences with the following
conditions:}
\end{theorem}

1) \textit{for each }$x\in X:=\tprod\limits_{i\in I}X_{i}$, clco$%
S_{i}(x)\subset T_{i}(x)$ and $S_{i}(x)\neq \emptyset $\textit{; }

2)\textit{\ }$S_{i}$\textit{\ is lower semicontinuous.}

\textit{Then there exists} $x^{\ast }\in D=\tprod\limits_{i\in I}D_{i}$ 
\textit{\ such that }$x_{i}^{\ast }\in T_{i}(x^{\ast })$ for each $i\in
I.\medskip $

In the present paper, our purpose is to give a fixed point theorem and to
research the equilibrium existence problem for abstract economies. In order
to establish our main results, we introduce the following definitions.

Let $X$ be a topological space, $Y$ be a nonempty subset of a topological
vector space $E$ and $D$ be a subset of $Y$.

\begin{definition}
The correspondence $T:X\rightarrow 2^{Y}$ is said to be \textit{w-lower
semicontinuous} (weakly lower semicontinuous) \textit{with respect to} $D$
if there exists a basis \ss\ of open symmetric neighbourhoods of $0$ in $E$
such that, for each $V\in $\ss , the correspondence $T^{V}$ is lower
semicontinuous, where $T^{V}(x)=(T(x)+V)\cap D$ for each $x\in X,$.
\end{definition}

\begin{remark}
By Lemma 2.6 in [19], it follows that if the correspondence $T:X\rightarrow
2^{Y}$ is almost lower semicontinuous, then it is w-lower semicontinuous
with respect to\textit{\ }$Y.$
\end{remark}

\begin{definition}
The correspondence $T:X\rightarrow 2^{Y}$ is said to be \textit{almost} 
\textit{w-lower semicontinuous} (almost weakly lower semicontinuous) \textit{%
with respect to} $D$ if there exists a basis \ss\ of open symmetric
neighbourhoods of $0$ in $E$ such that, for each $V\in $\ss , the
correspondence $\overline{T^{V}}$ is lower semicontinuous.\medskip
\end{definition}

\begin{example}
Let $\ T_{1}:(0,2)\rightarrow 2^{[1,4]}$ be the correspondence defined by
\end{example}

$T_{1}(x)=\left\{ 
\begin{array}{c}
\lbrack 2-x,2],\text{ if }x\in (0,1); \\ 
\{4\}\text{ \ \ \ \ \ \ if \ \ \ \ \ \ \ }x=1; \\ 
\lbrack 1,2]\text{ \ \ \ if \ \ \ }x\in (1,2).%
\end{array}%
\right. $

$T_{1}$ is not lower semicontinuous on $(0,2).$

Let $D=[1,2]$ and let $V=(-\varepsilon ,\varepsilon ),$ $\varepsilon >0,$ be
an open symmetric neighbourhood of $0$ in $\mathbb{R}.$ Then, it results that

for $\varepsilon \in (0,1),$

$T_{1}^{V}(x)=(T_{1}(x)+(-\varepsilon ,\varepsilon ))\cap D=\left\{ 
\begin{array}{c}
(2-x-\varepsilon ,2],\text{ if }x\in (0,1-\varepsilon ]; \\ 
\lbrack 1,2]\text{ \ \ if \ \ }x\in (1-\varepsilon ,1)\cup (1,2); \\ 
\phi \text{ \ \ \ \ \ \ \ \ \ \ \ \ \ \ \ \ if \ \ \ \ \ \ \ \ \ \ \ }x=1;%
\end{array}%
\right. ,$

for $\varepsilon \in \lbrack 1,2],$

$T_{1}^{V}(x)=(T_{1}(x)+(-\varepsilon ,\varepsilon ))\cap D=\left\{ 
\begin{array}{c}
\lbrack 1,2],\text{ if }x\in (0,1)\cup (1,2); \\ 
\phi \text{ \ \ \ \ \ \ \ \ \ if \ \ \ \ \ \ \ \ \ \ \ }x=1;%
\end{array}%
\right. $

if $\varepsilon \in (2,3],$

$T_{1}^{V}(x)=(T_{1}(x)+(-\varepsilon ,\varepsilon ))\cap D=\left\{ 
\begin{array}{c}
\lbrack 1,2]\text{\ \ if \ \ }x\in (0,1)\cup (1,2); \\ 
(4-\varepsilon ,2]\text{\ if \ \ \ \ \ \ \ \ \ \ \ \ \ \ \ }x=1;%
\end{array}%
\right. $

and if $\varepsilon >3,$

$T_{1}^{V}(x)=(T_{1}(x)+(-\varepsilon ,\varepsilon ))\cap D=[1,2]$\ if $x\in
(0,2).$

Then,

for $\varepsilon \in (0,1),$

$\overline{T_{1}^{V}}(x)=\left\{ 
\begin{array}{c}
\lbrack 2-x-\varepsilon ,2],\text{ if }x\in \lbrack 0,1-\varepsilon ); \\ 
\lbrack 1,2]\text{ \ \ \ \ if \ \ \ \ \ \ \ \ }x\in (1-\varepsilon ,2];%
\end{array}%
\right. $ and

for $\varepsilon \geq 1,$

$\overline{T_{1}^{V}}(x)=[1,2]$ for $x\in \lbrack 0,2].$

For each $V=(-\varepsilon ,\varepsilon )$ with $\varepsilon >0,$ the
correspondences $T_{1}^{V}$ and $\overline{T_{1}^{V}}$ are lower
semicontinuous and $\overline{T_{1}^{V}}$ has nonempty values. We conclude
that $T_{1}$ is w-lower semicontinuous with respect to $D$ and it is also
almost w-lower semicontinuous with respect to $D.$\medskip

\begin{proposition}
Let $X$ be a topological space, $Y$ be a nonempty subset of a topological
vector space $E.$ If the correspondence $T:X\rightarrow 2^{Y}$ is lower
semicontinuous and nonempty valued, then it is also w-lower semicontinuous
with respect to any set $D\subset Y$ with the property that $T(x)\cap D\neq
\phi ,$ for every $x\in X$.
\end{proposition}

\textit{Proof. }Let $V$ be an open symmetric neighborhood of $0$ in $E$.
Since the constant valued correspondence $x\rightarrow V$ is lower
semicontinuous, it follows that so it is the correspondence $x\rightarrow
(T(x)+V)$. Note also that this correspondence has nonempty open values and
that $(T(x)+V)\cap D\neq \emptyset $ for every $x\in X.$ Further,
Proposition 1 can be applied for $T_{1}(x)=T(x)+V$ and $T_{2}(x)=D$, $x\in
X.\medskip $

\begin{remark}
If the corespondence $T^{V}$ has empty values for some open set $V,$ it may
not be lower semicontinuous. The following example proves this assertion.
\end{remark}

\begin{example}
Let $T_{1}:[0,2]\rightarrow 2^{[3,5)}$ be the correspondence defined by
\end{example}

$T_{1}(x)=\left\{ 
\begin{array}{c}
\lbrack x+2,4]\text{ if }x\in \lbrack 0,1]; \\ 
(4,5)\text{ \ \ \ if \ \ \ }x\in (1,2];%
\end{array}%
\right. $.

Let $V=(-1,1)$ and $D=[0,3].$ Then,

$T_{1}^{(-1,1)}(x)=\left\{ 
\begin{array}{c}
(x+1,5)\text{ if }x\in \lbrack 0,1]; \\ 
(3,6)\text{ \ \ \ if \ \ \ }x\in (1,2];%
\end{array}%
\right. \cap $ $[0,3]=\left\{ 
\begin{array}{c}
(x+1,3]\text{ if }x\in \lbrack 0,1]; \\ 
\phi \text{ \ \ \ if \ \ \ }x\in (1,2];%
\end{array}%
\right. $

$T_{1}^{(-1,1)}$ is not lower semicontinuous.

\begin{remark}
$\overline{T^{V}}$ may not have convex values, even if $T^{V}$ is convex
valued.
\end{remark}

\begin{example}
Let $D=[1,2]$ and $T:[0,2]\rightarrow 2^{[0,4)}$ be the correspondence
defined by
\end{example}

$T(x)=\left\{ 
\begin{array}{c}
\lbrack 0,1]\text{ \ \ if }x\in \lbrack 0,1); \\ 
\phi \text{ \ \ \ \ \ \ \ \ if \ \ \ \ \ }x=1; \\ 
(2,3)\text{ \ \ \ if \ \ \ }x\in (1,2].%
\end{array}%
\right. $

$T$ is lower semicontinuous on [$0,2].$

If $\varepsilon \in (0,\frac{1}{2}),$ $T^{V}(x)=\left\{ 
\begin{array}{c}
\lbrack 1,1+\varepsilon )\text{ if }x\in \lbrack 0,1); \\ 
\phi \text{ \ \ \ \ \ \ if \ \ \ \ \ }x=1; \\ 
(2-\varepsilon ,2]\text{ \ \ \ if \ \ \ }x\in (1,2].%
\end{array}%
\right. $

Then, if $\varepsilon \in (0,\frac{1}{2}),$ $\overline{T^{V}}(x)=\left\{ 
\begin{array}{c}
\lbrack 1,1+\varepsilon ]\text{ \ \ \ \ \ \ \ \ \ \ \ if }x\in \lbrack 0,1);
\\ 
\lbrack 1,1+\varepsilon ]\cup \lbrack 2-\varepsilon ,2]\text{ if\ }x=1; \\ 
\lbrack 2-\varepsilon ,2]\text{ \ \ \ \ \ \ \ \ if \ \ \ }x\in (1,2];%
\end{array}%
\right. $ for $V=(-\varepsilon ,\varepsilon ).$

$\overline{T^{V}}$ does not have convex values in every point $x\in \lbrack
0,2].$

We also define the dual w-lower semicontinuity with respect to a set.

\begin{definition}
Let $T_{1},T_{2}:X\rightarrow 2^{Y}$ be correspondences. The pair $%
(T_{1},T_{2})$ is said to be \textit{dual almost w-lower semicontinuous}
(dual weakly lower semicontinuous) \textit{with respect to} $D$ if there
exists a basis \ss\ of open symmetric neighbourhoods of $0$ in $E$ such
that, for each $V\in $\ss , the correspondence $\overline{T_{(1,2)}^{V}}%
:X\rightarrow 2^{D}$ is lower semicontinuous, where $T_{(1,2)}^{V}:X%
\rightarrow 2^{D}$ is defined by $T_{(1,2)}^{V}(x)=(T_{1}(x)+V)\cap
T_{2}(x)\cap D$ for each $x\in X$.
\end{definition}

\begin{example}
Let $\ D=[1,2],$ $T_{1}:(0,2)\rightarrow 2^{[1,4]}$ be the correspondence
from the example 1 and $T_{2}:(0,2)\rightarrow 2^{[2,3]}$ be the
correspondence defined by
\end{example}

$T_{2}(x)=\left\{ 
\begin{array}{c}
\lbrack 2,3],\text{ if }x\in (0,1]; \\ 
\{2\}\text{ \ if \ \ }x\in (1,2);%
\end{array}%
\right. .$

The correspondences $T_{1}$ and $T_{2}$ are not semicontinuous.

For $\varepsilon \in (0,2],$ $(T_{1}(x)+(-\varepsilon ,\varepsilon ))\cap
D\cap T_{2}(x)=\left\{ 
\begin{array}{c}
\{2\}\text{ if }x\in (0,1)\cup (1,2); \\ 
\phi \text{ \ \ \ \ \ \ \ if \ \ \ \ \ \ \ \ \ }x=1.%
\end{array}%
\right. $

For $\varepsilon \in (2,\infty ),$ $(T_{1}(x)+(-\varepsilon ,\varepsilon
))\cap D\cap T_{2}(x)=\{2\}$ for each $x\in (0,2).$

Then, we have that for each $\varepsilon >0,$ $\overline{T_{(1,2)}^{V}}%
(x)=\{2\}$ for each $x\in \lbrack 0,2]$ and the correspondence $\overline{%
T_{(1,2)}^{V}}$ is lower semicontinuous and has nonempty values.

We conclude that the pair $(T_{1},T_{2})$ is dual almost w-lower
semicontinuous with respect to $D.$

\section{A new fixed point theorem}

We obtain a fixed point theorem which is an extension of Wu's fixed point
Theorem 1 in [18], in the sense that, for each $i\in I,$ the involved
correspondence $S_{i}$ is assumed to be almost w-lower semicontinuous with
respect to a set $D_{i},$ but $D_{i}$ is not convex as in the quoted result.

\begin{theorem}
\textit{Let }$I$\textit{\ be an index set. For each }$i\in I,$\textit{\ let }%
$X_{i}$\textit{\ be a nonempty convex subset of a Hausdorff locally convex
topological vector space }$E_{i}$\textit{, }$D_{i}$\ be a nonempty compact
convex metrizable subset of $X_{i}$ \textit{and }$S_{i},T_{i}:X:=\tprod%
\limits_{i\in I}X_{i}\rightarrow 2^{X_{i}}$\textit{\ be two correspondences
with the following conditions:}
\end{theorem}

1) \textit{for each }$x\in X$, $\overline{S}_{i}(x)\subset T_{i}(x)$\textit{%
. }

2)\textit{\ }$S_{i}$\textit{\ is almost w-lower semicontinuous with respect
to }$D_{i}$ \textit{and} $\overline{S_{i}^{V_{i}}}$ \textit{is} \textit{%
convex nonempty valued for each open absolutely convex symmetric
neighbourhood }$V_{i}$ \textit{of} $0$ \textit{in} $E_{i}$\textit{.}

\textit{Then there exists} $x^{\ast }\in D:=\tprod\limits_{i\in I}D_{i}$%
\textit{\ such that }$x_{i}^{\ast }\in T_{i}(x^{\ast })$ \textit{for each} $%
i\in I.\medskip $

\textit{Proof. }Since $D_{i}$ is compact, $D:=\tprod\limits_{i\in I}D_{i}$
is also compact in $X.$ For each $i\in I,$ let \ss $_{i}$\ be a basis of
open absolutely convex symmetric neighbourhoods of zero in $E_{i}$ and let 
\ss =$\tprod\limits_{i\in I}$\ss $_{i}.$ For each system of neighbourhoods $%
V=(V_{i})_{i\in I}\in \tprod\limits_{i\in I}$\ss $_{i},$ let's define the
corespondences $S_{i}^{V_{i}}:X\rightarrow 2^{D_{i}},$ by $%
S_{i}^{V_{i}}(x)=(S_{i}(x)+V_{i})\cap D_{i}$, $x\in X,$ $i\in I.$ By
assumption 2) each $\overline{S_{i}^{V_{i}}}$ is l.s.c with nonempty closed
convex values. According to Theorem 1.1 in [12], there exists a nonempty
valued, upper semicontinuous correspondence $G_{i}^{V_{i}}:D\rightarrow
2^{D_{i}}$ such that $G_{i}^{V_{i}}(x)\subset \overline{S_{i}^{V_{i}}}(x)$
for all $x\in D.$ Then, by Theorem 7.3.5 in [8] and Theorem 1.4 pag. 25 in
[21], the correspondence $F_{i}^{V_{i}}=$clco$G_{i}^{V_{i}}:D\rightarrow
2^{D_{i}}$ is also upper semicontinuous with nonempty closed convex values.
Let's define $F^{V}:D\rightarrow 2^{D}$ by $F^{V}(x)=\tprod\limits_{i\in
I}F_{i}^{V_{i}}(x)$ for each $x\in D.$ The correspondence $F^{V}$ is upper
semicontinuous with closed convex values. Therefore, according to
Himmelberg's fixed point theorem [7], there exists $x_{V}^{\ast
}=\tprod\limits_{i\in I}x_{V_{i}}^{\ast }\in D$ such that $x^{\ast }\in
F^{V}(x^{\ast }).$ It follows that $x_{V_{i}}^{\ast }\in \overline{%
S_{i}^{V_{i}}}(x_{V}^{\ast })$ for each $i\in I.$

For each $V=(V_{i})_{i\in I}\in $\ss $,$ let's define $Q_{V}=\cap _{i\in
I}\{x\in D:$ $x_{i}\in \overline{S_{i}^{V_{i}}}(x)\}.$

$Q_{V}$ is nonempty since $x_{V}^{\ast }\in Q_{V},$ then $Q_{V}$ is nonempty
and closed.

We prove that the family $\{Q_{V}:V\in \text{\ss }\}$ has the finite
intersection property.

Let $\{V^{(1)},V^{(2)},...,V^{(n)}\}$ be any finite set of $\text{\ss }$ and
let $V^{(k)}=\underset{i\in I}{\tprod }V_{i}^{(k)}$, $k=1,...,n.$ For each $%
i\in I$, let $V_{i}=\underset{k=1}{\overset{n}{\cap }}V_{i}^{(k)}$, then $%
V_{i}\in \text{\ss }_{i};$ thus $V=\underset{i\in I}{\tprod }V_{i}\in 
\underset{i\in I}{\tprod }\text{\ss }_{i}.$ Clearly $Q_{V}\subset \underset{%
k=1}{\overset{n}{\cap }}Q_{V^{(k)}}$ so that $\underset{k=1}{\overset{n}{%
\cap }}Q_{V^{(k)}}\neq \emptyset .$

Since $D$ is compact and the family $\{Q_{V}:V\in \text{\ss }\}$ has the
finite intersection property, we have that $\cap \{Q_{V}:V\in \text{\ss }%
\}\neq \emptyset .$ Take any $x^{\ast }\in \cap \{Q_{V}:V\in $\ss $\},$ then
for each $V_{i}\in \text{\ss }_{i},$ $x_{i}^{\ast }\in \overline{%
S_{i}^{V_{i}}}(x^{\ast })$. Acording to Lemma 1,\emph{\ }we have that\emph{\ 
} $x_{i}^{\ast }\in \overline{S_{i}}(x^{\ast }),$ for each $i\in I,$
therefore $x_{i}^{\ast }\in T(x^{\ast }).$\medskip

If $\left\vert I\right\vert =1$ we get the result bellow.

\begin{corollary}
\textit{Let }$X$\textit{\ be a nonempty convex subset of a Hausdorff locally
convex topological vector space }$F,$\textit{\ }$D$\textit{\ be a nonempty
compact convex metrizable subset of }$X$\textit{\ and }$S,T:X\rightarrow
2^{X}$\textit{\ be two correspondences with the following conditions:}
\end{corollary}

1) \textit{for each }$x\in X,$\textit{\ }$\overline{S}(x)\subset T(x)$%
\textit{\ and }$S(x)\neq \emptyset ,$

2)\textit{\ }$S$\textit{\ is almost w-lower semicontinuous with respect to }$%
D$ \textit{and }$\overline{S^{V}}$ \textit{is} \textit{convex nonempty
valued for each open absolutely convex symmetric neighbourhood }$V$ \textit{%
of} $0$ \textit{in} $E$.

\textit{Then, there exists a point }$x^{\ast }\in D$\textit{\ such that }$%
x^{\ast }\in T(x^{\ast }).$\medskip

In the particular case that the correspondence $S=T$ the following result
stands.

\begin{corollary}
Let $X$ be a nonempty convex subset of a \textit{Hausdorff locally convex
topological vector} space $F,$ $D$ be a nonempty compact convex metrizable
subset of $X$ and $T:X\rightarrow 2^{X}$ be an almost w- lower
semicontinuous correspondence with respect to $D$ and $\overline{T^{V}}$ 
\textit{is} \textit{convex nonempty valued for each open absolutely convex
symmetric neighborhood }$V$ \textit{of} $0$ \textit{in} $E$. Then, there
exists a point $x^{\ast }\in D$ such that $x^{\ast }\in \overline{T}(x^{\ast
}).\medskip $
\end{corollary}

\section{Applications in equilibrium theory}

Let $I$ be a nonempty set (the set of agents). For each $i\in I$, let $X_{i}$
be a nonempty topological vector space representing the set of actions and
define $X:=\underset{i\in I}{\prod }X_{i}$; let $A_{i}$, $B_{i}:X\rightarrow
2^{X_{i}}$ be the constraint correspondences and $P_{i}$ be the preference
correspondence for the agent $i$.

\begin{definition}[21]
\textit{\ }An \textit{abstract economy} $\Gamma
=(X_{i},A_{i},P_{i},B_{i})_{i\in I}$ is a family of ordered quadruples $%
(X_{i},A_{i},P_{i},B_{i})$.\medskip
\end{definition}

The notion of equilibrium plays a central role in the theory of equilibrium.
In the recent years the generalizations of this concept have been made in
some directions, several of them enlarging the set of acceptable points. One
of these methods is due to Yuan [21], who divided each constraint
correspondence into two parts, A and B, because the set of the fixed points
of \ the "small" correspondence may not be rich enough. Another method leads
us to the notion of "pseudo-equilibrium" and we will define it
further.\medskip

Here, for the definition of the equilibrium we follow Yuan [21].

\begin{definition}[21]
\textit{\ }An \textit{equilibrium} for $\Gamma $ is a point $x^{\ast }\in X$
such that for each $i\in I$, $x_{i}^{\ast }\in \overline{B}_{i}(x^{\ast })$
and $A_{i}(x^{\ast })\cap P_{i}(x^{\ast })=\emptyset $.
\end{definition}

\begin{remark}
When, for each $i\in I$, $A_{i}(x)=B_{i}(x)$ for all $x\in X,$ this abstract
economy model coincides with the classical one introduced by Borglin and
Keiding in [2]. If, in addition, $\overline{B}_{i}(x)=$cl$_{X_{i}}B_{i}(x)$
for each $x\in X,$ which is the case if $B_{i}$ has a closed graph in $%
X\times X_{i}$, the definition of equilibrium coincides with that one used
by Yannelis and Prabhakar in [20].
\end{remark}

\begin{remark}
An example of extension of the equilibrium model considering two constraint
correspondences $A_{i}$,$B_{i}:X\rightarrow 2^{X_{i}}$ for each player $i$
(with $A_{i}(x)\subset B_{i}(x)$, $x\in X$) is the notion of
quasi-equilibrium (see [6]) for an abstract economy, which has an analogue
in the private ownership economies. Even if Florenzano considers in [6] that
the interest of the quasi-equilibrium concept is purely mathematical, she
motivates the research of the conditions which guarantee its existence as
being very fruitful from a lot of points of view.
\end{remark}

The following example motivates the necessity of Yuan's model of abstract
economy with two constraint correspondences and illustrates it by using
correspondences for which the assumptions formulated by us hold, but those
made by Yannelis and Brahbakar [20] or by other authors (and which concern
the lower semicontinuity) do not hold.

\begin{example}
Let $\Gamma =(X,A,B,P)$ be an abstract economy with one agent, where $%
X=[0,4] $ and $A,B,P:X\rightarrow 2^{X}$ are defined below:
\end{example}

$A(x)=\left\{ 
\begin{array}{c}
\lbrack 1-x,2]\text{ if }x\in \lbrack 0,1); \\ 
\lbrack 3,4]\text{ \ \ \ \ \ if \ \ \ \ \ \ }x=1; \\ 
\lbrack 0,\frac{1}{2}]\text{ if\ \ \ \ \ \ \ }x\in (1,4]\text{;}%
\end{array}%
\right. $

$P(x)=\left\{ 
\begin{array}{c}
\lbrack \frac{3}{2},2]\text{ \ \ \ if \ }x\in \lbrack 0,1); \\ 
\{4\}\text{\ \ \ \ if \ \ \ \ \ \ \ \ \ \ }x=1; \\ 
\lbrack 1,2]\text{ \ \ \ if\ \ \ \ \ \ }x\in (1,4]\text{;}%
\end{array}%
\right. $

$B(x)=\left\{ 
\begin{array}{c}
\lbrack 1-x,2]\text{ if }x\in \lbrack 0,1); \\ 
\lbrack 3,4]\text{ \ \ if \ \ \ \ \ \ \ \ \ }x=1; \\ 
\lbrack 0,2]\text{ \ \ if\ \ \ \ \ \ }x\in (1,4]\text{.}%
\end{array}%
\right. $

The fixed point set of $A$ is Fix$(A)=[\frac{1}{2},1),$ the fixed point set
of $B$ is Fix$(B)=[\frac{1}{2},1)\cup (1,2]$ and $A(x)\subseteq B(x)$ for
each $x\in \lbrack 0,4].$ Since $U=\{x\in X:(A\cap P)(x)=\emptyset )=[1,4]$
and Fix$(A)\cap U=\emptyset ,$ Yannelis-Prahbakar's model ([20]) $(X,A,P)$
has not equilibrium points.

We notice that $x^{\ast }=\frac{3}{2}$ is an equilibrium point for Yuan's
model $(X,A,B,P):$ $(A\cap P)(\frac{3}{2})=\emptyset ,$ $\frac{3}{2}\in B(%
\frac{3}{2}).$ The correspondence $A$ proves to not have enough fixed points
and it must be enlarged by the correspondence $B.$

$B$, $P$ and $A$ are almost w-lower semicontinuous with respect to\textit{\ }%
$D=[0,2]$ and the game $(X,A,B,P)$ has equilibrium points, as we showed
above.

Since $P^{-1}(4)=\{x\in X:4\in P(x)\}=\{1\},$ the correspondence $P$ has not
open lower sections. We also see that the correspondences $A,B$ and $P$ are
not lower semicontinuous on $[0,4].$

There is large literature concerning the existence of the equilibrium in
Yuan's sense, which has been developed in the last decades. The authors
tried to generalize the properties of the involved correspondences; for an
overview, see for example [15] or [21]. In [21], Yuan provides applications
of abstract economies with two constraint correspondences to the systems of
generalized quasi-variational inequalities and to the systems of Ky Fan
minimax inequalities.

We define the following type of equilibrium for an abstract economy, which
is a slight extension of Yuan's equilibrium. The motivation of introducing
it is a mathematical one.

\begin{definition}
A \textit{pseudo} \textit{equilibrium} for $\Gamma $ is defined as a point $%
x^{\ast }\in X$ such that for each $i\in I$, $x_{i}^{\ast }\in \overline{B}%
_{i}(x^{\ast })$ and $x^{\ast }\in $cl$\{x\in X:(A_{i}\cap
P_{i})(x)=\emptyset \}$ for each $i\in I$.
\end{definition}

\begin{example}
Let $\Gamma $ $=(X,A,B,P)$ be an abstract economy with one agent, where $%
X=[0,4]$ and $A,B,P:X\rightarrow $ 2$^{X}$ are defined below:
\end{example}

$A(x)=\left\{ 
\begin{array}{c}
\lbrack 1-x,2]\text{ if }x\in \lbrack 0,1); \\ 
\lbrack 1,4]\text{ \ \ \ \ \ if \ \ \ \ \ \ }x=1; \\ 
\lbrack 0,\frac{1}{2}]\text{ if\ \ \ \ \ \ \ }x\in (1,4]\text{;}%
\end{array}%
\right. $

$P(x)=\left\{ 
\begin{array}{c}
\lbrack \frac{3}{2},2+x]\text{ if \ }x\in \lbrack 0,1); \\ 
\{1\}\text{\ \ \ \ if \ \ \ \ \ \ \ \ \ \ }x=1; \\ 
\lbrack 1,2]\text{ \ \ \ if\ \ \ \ \ \ }x\in (1,4]\text{;}%
\end{array}%
\right. $

$B(x)=\left\{ 
\begin{array}{c}
\lbrack 1-x,2]\text{ if }x\in \lbrack 0,1); \\ 
\lbrack 1,4]\text{ \ \ if \ \ \ \ \ \ \ \ \ }x=1; \\ 
\lbrack 0,2]\text{ \ \ if\ \ \ \ \ \ }x\in (1,4]\text{.}%
\end{array}%
\right. $

$(A\cap P)(x)=\left\{ 
\begin{array}{c}
\lbrack \frac{3}{2},2]\text{ \ \ if \ }x\in \lbrack 0,1); \\ 
\{1\}\text{ \ \ if \ \ \ \ \ \ \ \ \ }x=1; \\ 
\emptyset \text{ \ \ if\ \ \ \ \ \ }x\in (1,4]\text{.}%
\end{array}%
\right. $

We note that $x^{\ast }=1$ and $x^{\ast }=\frac{3}{2}$ are
pseudo-equilibrium points for $\Gamma ,$ since $1\in B(1),$ $\frac{3}{2}\in
B(\frac{3}{2})$ and $1,\frac{3}{2}\in $cl\{$x\in X:(A\cap P)(x)=\emptyset
\}=[1,4].\medskip $

Having a utility function $u_{i}:X\times X_{i}\rightarrow \mathbb{R}$ for
each agent $i$, we can define a preference correspondence $P_{i}$: $%
P_{i}(x):=\left\{ y_{i}\in X_{i}:u_{i}(x,y_{i})>u_{i}(x,x_{i})\right\} .$
Then, the condition of maximizing the utility function to obtain the
equilibrium point becomes: $A_{i}(x)\cap P_{i}(x)=\emptyset $ for each $i\in
I.\medskip $

Now we give an example of correspondence $P$ which is w-lower semicontinuous
with respect to a given set, $P$ being constructed from a utility function.

\begin{example}
Let $G$ be the game $([-1,1],A,P)$ with $I=\{1\}$, $A:[-1,1]\rightarrow
2^{[-1,1]}$ be defined as
\end{example}

$A(x)=\left\{ 
\begin{array}{c}
\lbrack 0,1],\text{ if }x\in \lbrack -1,0); \\ 
(\frac{1}{2},1]\text{ \ \ \ \ \ \ \ \ }x\in \lbrack 0,1].%
\end{array}%
\right. $ and

let $u:[-1,1]\times \lbrack -1,1]\rightarrow \mathbb{R}$ be the function
with levels defined as

$u(x,y)=\left\{ 
\begin{array}{c}
1,\text{\ \ \ \ if }(x,y)\in \lbrack -1,0)\times \lbrack -1,1)\cup \lbrack
0,1]\times \lbrack -1,0]\backslash \{(0,0)\}; \\ 
2\text{ \ \ \ \ \ \ \ \ \ \ \ \ \ \ \ \ \ if \ \ \ \ \ \ \ \ \ \ \ \ \ \ \ \
\ \ \ \ \ \ }(x,y)\in \lbrack -1,0)\times \{(1)\}; \\ 
3\text{ \ \ \ \ \ \ \ \ \ \ if \ \ \ \ \ \ \ \ \ \ \ \ \ \ \ \ \ \ \ \ }%
(x,y)\in \{(0,0)\}\cup \{0\}\times \lbrack 0,\frac{1}{2}); \\ 
4\text{ \ \ \ \ \ \ \ \ \ \ \ \ \ \ \ \ \ \ \ \ if \ \ \ \ \ \ \ \ \ \ \ \ \
\ \ \ \ \ \ \ \ \ \ \ }(x,y)\in \{0\}\times \lbrack \frac{1}{2},1]; \\ 
5\text{ \ \ \ \ \ \ \ \ \ \ \ \ \ \ \ \ \ \ \ \ \ \ if \ \ \ \ \ \ \ \ \ \ \
\ \ \ \ \ \ \ \ \ }(x,y)\in (0,1]\times (0,1); \\ 
6\text{ \ \ \ \ \ \ \ \ \ \ \ \ \ \ \ \ \ \ \ \ \ \ if \ \ \ \ \ \ \ \ \ \ \
\ \ \ \ \ \ \ \ \ \ \ }(x,y)\in (0,1]\times \{1\}.%
\end{array}%
\right. $

Then $P:[-1,1]\rightarrow 2^{[-1,1]}$ is defined as

$P(x):=\left\{ y\in X:u(x,y)>u(x,x)\right\} $

$P(x)=\left\{ 
\begin{array}{c}
\{1\},\text{ if }x\in \lbrack -1,0)\cup (0,1); \\ 
\lbrack \frac{1}{2},1]\text{ \ \ \ \ \ \ \ \ \ \ if \ \ \ \ \ \ \ }x=0; \\ 
\phi \text{ \ \ \ \ \ \ \ \ \ if \ \ \ \ \ \ \ \ \ \ \ }x=1.%
\end{array}%
\right. $

$P$ is not lower semicontinuous.

Let $V=(-\varepsilon ,\varepsilon )$ and $D=[1,2].$

$P(x)+V=\left\{ 
\begin{array}{c}
(1-\varepsilon ,1+\varepsilon ),\text{ if }x\in \lbrack -1,0)\cup (0,1); \\ 
(\frac{1}{2}-\varepsilon ,1+\varepsilon )\text{ \ \ \ \ \ \ \ \ \ \ \ if \ \
\ \ \ \ \ \ }x=0; \\ 
\phi \text{ \ \ \ \ \ \ \ \ \ \ \ \ \ \ \ \ \ \ if \ \ \ \ \ \ \ \ \ \ \ \ \ 
}x=1.%
\end{array}%
\right. $

For $\varepsilon \in (0,1],$

$P^{V}(x)=\left\{ 
\begin{array}{c}
\lbrack 1,1+\varepsilon ),\text{ if }x\in \lbrack -1,1); \\ 
\phi \text{ \ \ \ \ \ \ \ if \ \ \ \ \ \ \ \ \ }x=1.%
\end{array}%
\right. $

For $\varepsilon >1,$

$P^{V}(x)=\left\{ 
\begin{array}{c}
\lbrack 1,2],\text{ if }x\in \lbrack -1,1); \\ 
\phi \text{\ \ \ \ \ \ \ \ \ if\ \ \ \ }x=1.%
\end{array}%
\right. $

$P$ is w-lower semicontinuous with respect with $D=[1,2].$

$A(x)\cap P(x)=\left\{ 
\begin{array}{c}
\{1\},\text{ if }x\in \lbrack -1,0)\cup (0,1); \\ 
(\frac{1}{2},1]\text{ \ \ \ \ \ \ \ \ \ \ \ if \ \ \ \ \ \ }x=0; \\ 
\phi \text{ \ \ \ \ \ \ \ \ \ if \ \ \ \ \ \ \ \ \ \ }x=1.%
\end{array}%
\right. $

$A(1)\cap P(1)=\phi $, and, since $1\in A(1),$ we have that $x^{\ast }=1$ is
an equilibrium point for $G$.\medskip

As an application of the fixed point theorem proved in Section 3, we have
the following result.\medskip

\begin{theorem}
\textit{Let }$\Gamma =\left\{ X_{i},A_{i},B_{i},P_{i}\right\} _{i\in I}$%
\textit{\ be an abstract economy such that for each }$i\in I,$\textit{\ the
following conditions are fulfilled:}
\end{theorem}

1)\textit{\ }$X_{i}$\textit{\ is a nonempty convex subset of a Hausdorff
locally convex topological vector space }$E_{i}\,$\textit{\ and }$D_{i}$%
\textit{\ is a nonempty compact convex metrizable subset of }$X_{i}$\textit{;%
}

2)\textit{\ for each }$x\in X=\prod\limits_{i\in I}X_{i},$\textit{\ }$%
A_{i}\left( x\right) $\ \textit{and }$P_{i}(x)$\textit{\ are convex} \textit{%
and }$B_{i}\left( x\right) $\textit{\ is nonempty, convex and }$A_{i}\left(
x\right) \subset B_{i}(x)$\textit{;}

3) $W_{i}=\left\{ x\in X:A_{i}\left( x\right) \cap P_{i}\left( x\right) \neq
\emptyset \right\} $\textit{\ is closed in }$X$\textit{;}

4)\textit{\ }$H_{i}:X\rightarrow 2^{X_{i}}$\textit{\ defined by }$%
H_{i}\left( x\right) =A_{i}(x)\cap P_{i}\left( x\right) $\textit{\ for each }%
$x\in X$\textit{\ is almost w-lower semicontinuous with respect to }$D_{i}$ 
\textit{on }$W_{i}$ and $\overline{H_{i}^{V_{i}}}$ \textit{is} \textit{%
convex nonempty valued for each open absolutely convex symmetric
neighbourhood }$V_{i}$ \textit{of} $0$ \textit{in} $E_{i}$\textit{;}

5)\textit{\ }$B_{i}:X\rightarrow 2^{X_{i}}$\textit{\ is almost w-lower
semicontinuous with respect to }$D_{i}$ and $\overline{B_{i}^{V_{i}}}$ 
\textit{is} \textit{convex nonempty valued for each open absolutely convex
symmetric neighbourhood }$V_{i}$ \textit{of} $0$ \textit{in} $E_{i}$\textit{;%
}

6)\textit{\ for each }$x\in X$\textit{\ , }$x_{i}\notin \overline{\mathit{(}%
A_{i}\cap P_{i})}\left( x\right) $\textit{.}

\textit{Then there exists }$x^{\ast }\in D=$\textit{\ }$\prod\limits_{i\in
I}D_{i}$\textit{\ such that }$x_{i}^{\ast }\in \overline{B}_{i}\left(
x^{\ast }\right) $\textit{\ and }$x^{\ast }\in $cl$\{x\in X:(A_{i}\cap
P_{i})(x)=\emptyset \}$ for each $i\in I$ ($x^{\ast }$ is a pseudo
equilibrium point for $\Gamma $)$.\medskip $

\textit{Proof. }Let $i\in I.$ \ By condition 3) we know that $W_{i}$ is
closed in $X.$

Let's define $T_{i}:X\rightarrow 2^{X_{i}}$ by $T_{i}\left( x\right)
=\left\{ 
\begin{array}{c}
A_{i}\left( x\right) \cap P_{i}\left( x\right) ,\text{ if }x\in W_{i}, \\ 
B_{i}\left( x\right) ,\text{ \ \ \ \ \ \ \ \ \ \ \ if }x\notin W_{i}%
\end{array}%
\right. $ for each $x\in X.$

Then $T_{i}:X\rightarrow 2^{X_{i}}$ is a correspondence with nonempty convex
values. We shall prove that $T_{i}:X\rightarrow 2^{D_{i}}$ is almost w-lower
semicontinuous with respect to $D_{i}$. Let \ss $_{i}$\ be a basis of open
absolutely convex symmetric neighbourhoods of $0$ in $E_{i}$ and let \ss =$%
\tprod\limits_{i\in I}$\ss $_{i}.$ For each $V=(V_{i})_{i\in I}\in
\tprod\limits_{i\in I}$\ss $_{i},$ for each $x\in X,$ let for each $i\in I$

$B^{V_{i}}(x)=(B_{i}\left( x\right) +V_{i})\cap D_{i}$,

$F^{V_{i}}(x)=((A_{i}\left( x\right) \cap P_{i}\left( x\right) )+V_{i})\cap
D_{i}$ and

$T_{i}^{V_{i}}(x)=\left\{ 
\begin{array}{c}
F^{V_{i}}(x),\text{ if }x\in W_{i}, \\ 
B^{V_{i}}(x),\text{\ if }x\notin W_{i}.%
\end{array}%
\right. $

For each closed set $V_{i}^{\prime }$ in $D_{i}$, the set

$\left\{ x\in X:\overline{T_{i}^{V_{i}}}\left( x\right) \subset
V_{i}^{\prime }\right\} =$

$=\left\{ x\in W_{i}:\overline{F^{V_{i}}}(x)\subset V_{i}^{\prime }\right\}
\cup \left\{ x\in X\smallsetminus W_{i}:\overline{B^{V_{i}}}(x)\subset
V_{i}^{^{\prime }}\right\} $

$=\left\{ x\in W_{i}:\overline{F^{V_{i}}}(x)\subset V_{i}^{^{\prime
}}\right\} \cup \left\{ x\in X:\overline{B^{V_{i}}}(x)\subset V_{i}^{\prime
}\right\} .$

According to condition 6), the set $\left\{ x\in W_{i}:\overline{F^{V_{i}}}%
(x)\subset V_{i}^{\prime }\right\} $ is closed in $X$. The set $\left\{ x\in
X:\overline{B^{V_{i}}}(x)\subset V_{i}^{\prime }\right\} $ is closed in $X$
because $\overline{B^{V_{i}}}$ is lower semicontinuous.

Therefore, the set $\left\{ x\in X:\overline{T_{i}^{V_{i}}}\left( x\right)
\subset V_{i}^{\prime }\right\} $ is closed in $X.$ It shows that $\overline{%
T_{i}^{V_{i}}}:X\rightarrow 2^{D_{i}}$ is lower semicontinuous. According to
Theorem 2, there exists $x^{\ast }\in D=$ $\prod\limits_{i\in I}D_{i}$ such
that $x^{\ast }\in \overline{T}_{i}\left( x^{\ast }\right) ,$ for each $i\in
I.$ By condition 6) we have that $x_{i}^{\ast }\in \overline{B}_{i}\left(
x^{\ast }\right) $ and $x^{\ast }\in $cl$\{x\in X:(A_{i}\cap
P_{i})(x)=\emptyset \}$ for each $i\in I.\medskip $

\begin{example}
Let $\Gamma =\left\{ X_{i},A_{i},B_{i},P_{i}\right\} _{i\in I}$\textit{\ }be
an abstract economy, where $I=\{1,2,...,n\},$ $X_{i}=[0;4]$ be a compact
convex choice set, $D_{i}=[0,2]$ for each $i\in I$ and $X=\prod\limits_{i\in
I}X_{i}$.
\end{example}

Let $A=\{x\in X:$ $\forall j\in \{1,2,...,n\},$ $x_{j}\in \lbrack 0,1]$ and $%
\exists j\in \{1,2,...,n\}$ such that $x_{j}=1\}$ and let the
correspondences $A_{i},B_{i},P_{i}:X\rightarrow 2^{X_{i}}$ be defined as
follows:

for each $(x_{1},x_{2},...,x_{n})\in X,$

$A_{i}(x)=\left\{ 
\begin{array}{c}
\lbrack 1-x_{i},2]\text{ if }x\in X\text{ and }\forall j\in \{1,2,...,n\},%
\text{ }x_{j}\in \lbrack 0,1); \\ 
\lbrack 3,4]\text{ \ \ \ \ \ \ \ \ \ \ \ \ \ \ \ \ \ \ \ \ \ \ \ \ \ if \ \
\ \ \ \ \ \ \ \ \ \ \ \ \ \ \ \ \ \ \ \ \ \ \ \ }x\in A; \\ 
\lbrack 0,\frac{1}{2}]\text{, \ \ \ \ \ \ \ \ \ \ \ \ \ \ \ \ \ \ \ \ \ \ \
\ \ \ \ \ \ \ \ \ \ \ \ \ \ \ \ \ \ \ \ \ \ \ \ otherwise;}%
\end{array}%
\right. $

$P_{i}(x)=\left\{ 
\begin{array}{c}
\lbrack \frac{3}{2},2+x_{i}]\text{ if }x\in X\text{ and }\forall j\in
\{1,2,...,n\},\text{ }x_{j}\in \lbrack 0,1); \\ 
\{4\}\text{ \ \ \ \ \ \ \ \ \ \ \ \ \ \ \ \ \ \ \ \ \ if \ \ \ \ \ \ \ \ \ \
\ \ \ \ \ \ \ \ \ \ \ \ \ \ \ \ \ \ \ \ }x\in A; \\ 
\lbrack 1,2]\text{, \ \ \ \ \ \ \ \ \ \ \ \ \ \ \ \ \ \ \ \ \ \ \ \ \ \ \ \
\ \ \ \ \ \ \ \ \ \ \ \ \ \ \ \ \ \ otherwise;}%
\end{array}%
\right. $

$B_{i}(x)=\left\{ 
\begin{array}{c}
\lbrack 1-x_{i},2]\text{ if }x\in X\text{ and }\forall j\in \{1,2,...,n\},%
\text{ }x_{j}\in \lbrack 0,1); \\ 
\lbrack 3,4]\text{ \ \ \ \ \ \ \ \ \ \ \ \ \ \ \ \ \ \ \ \ \ \ if \ \ \ \ \
\ \ \ \ \ \ \ \ \ \ \ \ \ \ \ \ \ \ \ \ \ \ \ }x\in A; \\ 
\lbrack 0,2]\text{, \ \ \ \ \ \ \ \ \ \ \ \ \ \ \ \ \ \ \ \ \ \ \ \ \ \ \ \
\ \ \ \ \ \ \ \ \ \ \ \ \ \ \ \ \ \ \ otherwise.}%
\end{array}%
\right. $

The correspondences $A_{i},B_{i},P_{i}$ are not lower semicontinuous on $X.$

$A_{i}(x)\cap P_{i}(x)=\left\{ 
\begin{array}{c}
\lbrack \frac{3}{2},2]\text{ if }x\in X\text{ and }\forall j\in
\{1,2,...,n\},\text{ }x_{j}\in \lbrack 0,1); \\ 
\{4\}\text{ \ \ \ \ \ \ \ \ \ \ \ \ \ \ \ \ \ if \ \ \ \ \ \ \ \ \ \ \ \ \ \
\ \ \ \ \ \ \ \ \ \ \ \ \ \ }x\in A; \\ 
\phi \text{, \ \ \ \ \ \ \ \ \ \ \ \ \ \ \ \ \ \ \ \ \ \ \ \ \ \ \ \ \ \ \ \
\ \ \ \ \ \ \ \ \ \ \ otherwise.}%
\end{array}%
\right. $

$W_{i}=\left\{ x\in X:A_{i}\left( x\right) \cap P_{i}\left( x\right) \neq
\emptyset \right\} =[0,1]^{n}$\textit{\ }is closed in\textit{\ }$X.$

$\overline{\mathit{(}A_{i}\cap P_{i})}\left( x\right) =\left\{ 
\begin{array}{c}
\lbrack \frac{3}{2},2]\text{ if }x\in X\text{ and }\forall j\in
\{1,2,...,n\},\text{ }x_{j}\in \lbrack 0,1); \\ 
\lbrack \frac{3}{2},2]\cup \{4\}\text{ \ \ \ \ \ \ \ \ \ \ \ \ \ \ \ \ \ if\
\ \ \ \ \ \ \ \ \ \ \ \ \ \ \ \ \ \ }x\in A; \\ 
\phi \text{, \ \ \ \ \ \ \ \ \ \ \ \ \ \ \ \ \ \ \ \ \ \ \ \ \ \ \ \ \ \ \ \
\ \ \ \ \ \ \ \ \ \ \ \ otherwise.}%
\end{array}%
\right. $

We notice that for each\textit{\ }$x\in X$\textit{\ , }$x_{i}\notin 
\overline{\mathit{(}A_{i}\cap P_{i})}\left( x\right) .$

We shall prove that $B_{i}$\ and $\mathit{(}A_{i}\cap P_{i})_{W_{i}}$ are
almost w-lower semicontinuous with respect to\textit{\ }$D_{i}=[0,2].$

On $W_{i},$

$\mathit{(}A_{i}\cap P_{i})\left( x\right) =\left\{ 
\begin{array}{c}
\lbrack \frac{3}{2},2]\text{ if }x\in X\text{, }\forall j\in \{1,2,...,n\},%
\text{ }x_{j}\in \lbrack 0,1); \\ 
\{4\}\text{ \ \ \ \ \ \ \ \ \ \ \ \ \ \ \ \ \ \ \ \ if \ \ \ \ \ \ \ \ \ \ \
\ \ \ \ \ \ \ \ }x\in A;%
\end{array}%
\right. $

$\mathit{(}A_{i}\cap P_{i})\left( x\right) +(-\varepsilon ,\varepsilon
)=\left\{ 
\begin{array}{c}
(\frac{3}{2}-\varepsilon ,2+\varepsilon )\text{ if }x\in X\text{, }\forall
j\in \{1,2,...,n\},\text{ }x_{j}\in \lbrack 0,1); \\ 
(4-\varepsilon ,4+\varepsilon )\text{ \ \ \ \ \ \ \ \ \ \ \ \ \ \ \ \ \ if \
\ \ \ \ \ \ \ \ \ \ \ \ \ \ \ \ \ \ }x\in A;%
\end{array}%
\right. $

Let $\mathit{(}A_{i}\cap P_{i})^{V}\left( x\right) =(\mathit{(}A_{i}\cap
P_{i})\left( x\right) +(-\varepsilon ,\varepsilon ))\cap \lbrack 0,2],$
where $V=(-\varepsilon ,\varepsilon ).$

Then,

if $\varepsilon \in (0,\frac{3}{2}],$

$\mathit{(}A_{i}\cap P_{i})^{V}\left( x\right) =\left\{ 
\begin{array}{c}
(\frac{3}{2}-\varepsilon ,2]\text{ if }x\in X\text{ and }\forall j\in
\{1,2,...,n\},\text{ }x_{j}\in \lbrack 0,1); \\ 
\phi \text{ \ \ \ \ \ \ \ \ \ \ \ \ \ \ \ \ \ \ \ \ \ \ \ if \ \ \ \ \ \ \ \
\ \ \ \ \ \ \ \ \ \ \ \ \ \ \ \ \ \ \ \ }x\in A;%
\end{array}%
\right. $

if $\varepsilon \in (\frac{3}{2},2],$

$\mathit{(}A_{i}\cap P_{i})^{V}\left( x\right) =\left\{ 
\begin{array}{c}
\lbrack 0,2]\text{ if }x\in X\text{ and }\forall j\in \{1,2,...,n\},\text{ }%
x_{j}\in \lbrack 0,1); \\ 
\phi \text{ \ \ \ \ \ \ \ \ \ \ \ \ \ \ \ \ \ \ \ \ \ if \ \ \ \ \ \ \ \ \ \
\ \ \ \ \ \ \ \ \ \ \ \ \ \ \ \ }x\in A;%
\end{array}%
\right. $

if $\varepsilon \in (2,4],$

$\mathit{(}A_{i}\cap P_{i})^{V}\left( x\right) =\left\{ 
\begin{array}{c}
\lbrack 0,2]\text{ if }x\in X\text{ and }\forall j\in \{1,2,...,n\},\text{ }%
x_{j}\in \lbrack 0,1); \\ 
\text{ }(4-\varepsilon ,2]\text{ \ \ \ \ \ \ \ \ \ \ \ \ \ if \ \ \ \ \ \ \
\ \ \ \ \ \ \ \ \ \ \ \ \ \ \ \ \ \ }x\in A;%
\end{array}%
\right. $

and if $\varepsilon >4,$

$\mathit{(}A_{i}\cap P_{i})^{V}\left( x\right) =[0,2]$ if $x\in \lbrack
0,1]^{n}.$

Hence, for each $V=(-\varepsilon ,\varepsilon ),$ $\overline{\mathit{(}%
A_{i}\cap P_{i})^{V}}_{W_{i}}$ is lower semicontinuous and has nonempty
values.

$B_{i}\left( x\right) +(-\varepsilon ,\varepsilon )=\left\{ 
\begin{array}{c}
(1-x_{i}-\varepsilon ,\text{ }2+\varepsilon )\text{ if }x\in X\text{, }%
\forall j\in \{1,2,...,n\},\text{ }x_{j}\in \lbrack 0,1); \\ 
(3-\varepsilon ,\text{ }4+\varepsilon )\text{ \ \ \ \ \ \ \ \ \ \ \ \ \ \ \
\ \ \ \ \ \ \ \ \ \ \ \ if \ \ \ \ \ \ \ \ \ \ \ \ \ \ \ \ }x\in A; \\ 
(-\varepsilon ,\text{ }2+\varepsilon )\text{ \ \ \ \ \ \ \ \ \ \ \ \ \ \ \ \
\ \ \ \ \ \ \ \ \ \ \ \ \ \ \ \ \ \ \ \ \ \ \ \ \ \ \ \ otherwise.}%
\end{array}%
\right. $

Let $B_{i}{}^{V}\left( x\right) =(B_{i}\left( x\right) +(-\varepsilon
,\varepsilon ))\cap \lbrack 0,2],$ where $V=(-\varepsilon ,\varepsilon ).$

Then,

if $\varepsilon \in (0,1],$

$B_{i}{}^{V}\left( x\right) =\left\{ 
\begin{array}{c}
(1-x_{i}-\varepsilon ,2]\text{ \ \ \ \ \ \ \ \ \ \ \ \ \ \ \ \ \ \ \ \ \ \ \
\ \ \ if }x\in X,\text{ }x_{i}\in \lbrack 0,1-\varepsilon ]\text{, } \\ 
\text{ \ \ \ \ \ \ \ \ \ \ \ \ \ \ \ \ \ \ \ \ \ \ \ \ \ \ \ \ \ \ }\forall
j\in \{1,2,...,n\}\setminus \{i\},\text{ }x_{j}\in \lbrack 0,1); \\ 
\lbrack 0,2]\text{ \ \ \ \ \ \ \ \ \ \ \ \ \ \ \ \ \ \ \ \ \ \ \ \ \ \ \ \ \
\ \ \ \ \ \ \ \ \ if }x\in X,\text{ }x_{i}\in (1-\varepsilon ,1)\text{,} \\ 
\text{ \ \ \ \ \ \ \ \ \ \ \ \ \ \ \ \ \ \ \ \ \ \ \ \ \ \ \ \ \ \ \ \ }%
\forall j\in \{1,2,...,n\}\setminus \{i\},\text{ }x_{j}\in \lbrack 0,1); \\ 
\phi \text{ \ \ \ \ \ \ \ \ \ \ \ \ \ \ \ \ \ \ \ \ \ \ \ \ \ \ \ \ \ \ \ \
\ \ \ \ \ \ \ \ \ \ \ \ \ \ \ \ \ if \ \ \ \ \ \ \ \ \ \ \ \ \ \ \ }x\in A;
\\ 
\lbrack 0,2]\text{ \ \ \ \ \ \ \ \ \ \ \ \ \ \ \ \ \ \ \ \ \ \ \ \ \ \ \ \ \
\ \ \ \ \ \ \ \ \ \ \ \ \ \ \ \ \ \ \ \ \ \ \ \ \ \ \ \ \ otherwise;}%
\end{array}%
\right. $

if $\varepsilon \in (1,3],$ $B_{i}{}^{V}\left( x\right) =\left\{ 
\begin{array}{c}
\lbrack 0,2]\text{ if }x\in X\text{ and }\forall j\in \{1,2,...,n\},\text{ }%
x_{j}\in \lbrack 0,1); \\ 
(3-\varepsilon ,2]\text{ \ \ \ \ \ \ \ \ \ \ \ \ \ \ \ \ \ if \ \ \ \ \ \ \
\ \ \ \ \ \ \ \ \ \ \ \ \ \ }x\in A; \\ 
\lbrack 0,2]\text{ \ \ \ \ \ \ \ \ \ \ \ \ \ \ \ \ \ \ \ \ \ \ \ \ \ \ \ \ \
\ \ \ \ \ \ \ \ \ \ \ \ otherwise.}%
\end{array}%
\right. $

and if $\varepsilon >3,$ $B_{i}{}^{V}\left( x\right) =[0,2]$ if $x\in X.$

Then, for each $V=(-\varepsilon ,\varepsilon ),$ $\overline{B_{i}^{V}}$ is
lower semicontinuous and has nonempty values.

Therefore, all hypotheses of Theorem 3 are satisfied, so that there exists
an equilibrium point $x^{\ast }=\{\frac{3}{2},\frac{3}{2},...,\frac{3}{2}%
\}\in X$ such that $x_{i}^{\ast }\in \overline{B}_{i}\left( x^{\ast }\right) 
$ and $(A_{i}\cap P_{i})(x^{\ast })=\emptyset .$\bigskip

Theorem 4 deals with abstract economies which have dual w-lower
semicontinuous pairs of correspondences and can be compared with Theorem 5
in Wu [19].

\begin{theorem}
\textit{Let }$\Gamma =\left\{ X_{i},A_{i},B_{i},P_{i}\right\} _{i\in I}$%
\textit{\ be an abstract economy such that for each }$i\in I,$\textit{\ the
following conditions are fulfilled:}
\end{theorem}

1)\textit{\ }$X_{i}$\textit{\ is a nonempty convex subset of a Hausdorff
locally convex topological vector space }$E_{i}\,$and\textit{\ }$D_{i}$%
\textit{\ is a nonempty compact convex metrizable subset of }$X_{i}$\textit{;%
}

2)\textit{\ for each }$x\in X=\prod\limits_{i\in I}X_{i},$\textit{\ }$%
P_{i}(x)\subset D_{i},$\textit{\ and }$B_{i}\left( x\right) $\textit{\ is
nonempty;}

3)\textit{\ the set }$W_{i}=\left\{ x\in X:A_{i}\left( x\right) \cap
P_{i}\left( x\right) \neq \emptyset \right\} $\textit{\ is open in }$X$%
\textit{;}

4)\textit{\ the pair }$(A_{i\mid \text{cl}W_{i}},P_{i\mid \text{cl}W_{i}})$ 
\textit{is dual almost w-lower semicontinuous with respect to }$D_{i}$%
\textit{, }$B_{i}:X\rightarrow 2^{X_{i}}$\textit{\ is almost w-lower
semicontinuous with respect to }$D_{i}$;

5) \textit{if} $T_{i,V_{i}}:X\rightarrow 2^{X_{i}}$ \textit{is defined by} $%
T_{i,V_{i}}(x)=(A_{i}(x)+V_{i})\cap D_{i}\cap P_{i}(x)$ \textit{for each} $%
x\in X,$ \textit{then the} \textit{correspondences} $\overline{B_{i}^{V_{i}}}
$ \textit{and} $\overline{T_{i,V_{i}}}$ \textit{are} \textit{nonempty} 
\textit{convex valued for each open absolutely convex symmetric
neighbourhood }$V_{i}$ \textit{of} $0$ \textit{in} $E_{i}$\textit{;}

6)\textit{\ for each }$x\in X$\textit{\ , }$x_{i}\notin \overline{P}%
_{i}\left( x\right) $\textit{;}

\textit{Then, there exists }$x^{\ast }\in D=$\textit{\ }$\prod\limits_{i\in
I}D_{i}$\textit{\ such that }$x_{i}^{\ast }\in \overline{B}_{i}\left(
x^{\ast }\right) $\textit{\ and }$A_{i}\left( x^{\ast }\right) \cap
P_{i}\left( x^{\ast }\right) =\emptyset $\textit{\ for all }$i\in I.\medskip 
$

\textit{Proof.} For each $i\in I,$ let \ss $_{i}$\ denote the family of all
open absolutely convex symmetric neighbourhoods of zero in $E_{i}$ and let 
\ss $=\tprod\limits_{i\in I}$\ss $_{i}.$ For each $V=\tprod\limits_{i\in
I}V_{i}\in \tprod\limits_{i\in I}$\ss $_{i},$ for each $i\in I,$ let

$B^{V_{i}}(x)=(B_{i}\left( x\right) +V_{i})\cap D_{i}$ for each $x\in X$ and

$S_{i}^{V_{i}}\left( x\right) =\left\{ 
\begin{array}{c}
T_{i,V_{i}}(x),\text{ \ \ \ \ \ \ \ \ \ \ \ \ \ \ \ \ \ \ if }x\in \text{cl}%
W_{i}, \\ 
B_{i}^{V_{i}}(x),\text{ \ \ \ \ \ \ \ \ \ \ \ \ \ \ \ \ if }x\notin \text{cl}%
W_{i},%
\end{array}%
\right. $

$\overline{S_{i}^{V_{i}}}$ has closed values. Next, we shall prove that $%
\overline{S_{i}^{V_{i}}}:X\rightarrow 2^{D_{i}}$ is lower semicontinuous.

For each closed set $V^{\prime }$ in $D_{i}$, the set

$\left\{ x\in X:\overline{S_{i}^{V_{i}}}\left( x\right) \subset V^{\prime
}\right\} =$

$=\left\{ x\in \text{cl}W_{i}:\overline{T_{i,V_{i}}}(x)\subset V^{\prime
}\right\} \cup \left\{ x\in X\smallsetminus \text{cl}W_{i}:\overline{%
B_{i}^{V_{i}}}(x)\subset V^{\prime }\right\} $

=$\left\{ x\in \text{cl}W_{i}:\overline{T_{i,V_{i}}}(x)\subset V^{\prime
}\right\} \cup \left\{ x\in X:\overline{B_{i}^{V_{i}}}(x)\subset V^{\prime
}\right\} .$

We know that the correspondence $\overline{T_{i,V_{i}}}(x)_{\mid \text{cl}%
W_{i}}:$ cl$W_{i}\rightarrow 2^{D_{i}}$ is lower semicontinuous. The set $%
\left\{ x\in \text{cl}W_{i}:\overline{T_{i,V_{i}}}(x)\subset V^{\prime
}\right\} $ is \ closed in $\ $cl$W_{i}$, and hence it is also closed in $X$
because cl$W_{i}$ is closed in $X$. Since $\overline{B_{i}^{V_{i}}}%
(x):X\rightarrow 2^{D_{i}}$ is lower semicontinuous, the set $\{x\in X:%
\overline{B_{i}^{V_{i}}}(x)\}\subset V^{\prime }$ is closed in $X$ and
therefore the set $\left\{ x\in X:\overline{S_{i}^{V_{i}}}\left( x\right)
\subset V^{\prime }\right\} $ is closed in $X$. It showes that $\overline{%
S_{i}^{V_{i}}}:X\rightarrow 2^{D_{i}}$ is lower semicontinuous. According to
Wu's Theorem 1, applied for the correspondences $\overline{S_{i}^{V_{i}}}%
=T_{i}^{V_{i}},$ there exists a point $x_{V}^{\ast }\in D=$ $%
\prod\limits_{i\in I}D_{i}$ such that $x_{V_{i}}^{\ast }\in
T_{i}^{V_{i}}\left( x_{V}^{\ast }\right) $ for each $i\in I.$ By condition
5), we have that $x_{V_{i}}^{\ast }\notin \overline{P_{i}}\left( x_{V}^{\ast
}\right) ,$ hence, $x_{Vi}^{\ast }\notin \overline{A_{i}^{V_{i}}}\left(
x_{V}^{\ast }\right) \cap \overline{P_{i}}\left( x_{V}^{\ast }\right) $. 
\newline
We also have that clGr$(T_{i,V_{i}})\subseteq $ clGr$(A_{i}^{V_{i}})\cap $%
clGr$P_{i}.$Then $\overline{T_{i,V_{i}}}(x)\subseteq $ $\overline{%
A_{i}^{V_{i}}}(x)\cap \overline{P_{i}}\left( x\right) $ for each $x\in X.$
It follows that $x_{Vi}^{\ast }\notin \overline{T_{i,V_{i}}}(x_{V}^{\ast }).$
Therefore, $x_{Vi}^{\ast }\in \overline{B^{V_{i}}}\left( x_{V}^{\ast
}\right) .$

For each $V=(V_{i})_{i\in I}\in \tprod\limits_{i\in I}$\ss $_{i},$ let's
define $Q_{V}=\cap _{i\in I}\{x\in D:x\in \overline{B^{V_{i}}}\left(
x\right) $ and $A_{i}\left( x\right) \cap P_{i}\left( x\right) =\emptyset
\}. $

$Q_{V}$ is nonempty since $x_{V}^{\ast }\in Q_{V},$ and it is a closed
subset of $D$ according to 3). Then, $Q_{V}$ is nonempty and compact.

Let \ss =$\tprod\limits_{i\in I}$\ss $_{i}.$ We prove that the family $%
\{Q_{V}:V\in \text{\ss }\}$ has the finite intersection property.

Let $\{V^{(1)},V^{(2)},...,V^{(n)}\}$ be any finite set of $\text{\ss\ }$and
let $V^{(k)}=\underset{i\in I}{\tprod }V_{i}^{(k)}{}_{i\in I}$, $k=1,...,n.$
For each $i\in I$, let $V_{i}=\underset{k=1}{\overset{n}{\cap }}V_{i}^{(k)}$%
, then $V_{i}\in \text{\ss }_{i};$ thus $V\in \underset{i\in I}{\tprod }%
\text{\ss }_{i}.$ Clearly $Q_{V}\subset \underset{k=1}{\overset{n}{\cap }}%
Q_{V^{(k)}}$ so that $\underset{k=1}{\overset{n}{\cap }}Q_{V^{(k)}}\neq
\emptyset .$

Since $D$ is compact and the family $\{Q_{V}:V\in \text{\ss }\}$ has the
finite intersection property, we have that $\cap \{Q_{V}:V\in \text{\ss }%
\}\neq \emptyset .$ Take any $x^{\ast }\in \cap \{Q_{V}:V\in $\ss $\},$ then
for each $V\in \text{\ss },$

$x^{\ast }\in \cap _{i\in I}\left\{ x^{\ast }\in D:x_{i}^{\ast }\in 
\overline{B^{V_{i}}}\left( x\right) \text{ and }A_{i}\left( x\right) \cap
P_{i}\left( x\right) =\emptyset )\right\} .$

Hence, $x_{i}^{\ast }\in \overline{B^{V_{i}}}\left( x^{\ast }\right) $ for
each $V\in $\ss\ and for each $i\in I.$ According to Lemma\emph{\ }1,\emph{\ 
}we have that\emph{\ } $x_{i}^{\ast }\in \overline{B_{i}}(x^{\ast })$ and $%
(A_{i}\cap P_{i})(x^{\ast })=\emptyset $ for each $i\in I.$\medskip

Now we introduce the next concept which also generalizes the lower
semicontinuous correspondences.

Let $X$ be a non-empty convex subset of a topological linear space $E$, $Y$
be a non-empty set in a topological space and $K\subseteq X\times Y.$

\begin{definition}
The correspondence $T:X\times Y\rightarrow 2^{X}$ has the e-LSCS-property
(e-lower semicontinuous selection property) on $K,$ if for each absolutely
convex neighborhood $V$ of zero in $E,$ there exists a lower semicontinuous
correspondence with convex values $S^{V}:X\times Y\rightarrow 2^{X}$ such
that $S^{V}(x,y)\subset T(x,y)+$cl$V$ and $x\notin $cl$S^{V}(x,y)$ for every 
$(x,y)\in K$.$\medskip $
\end{definition}

The following theorem concerns the abstract economies which have
correspondences with e-LSCS-property.

\begin{theorem}
\textit{Let }$\Gamma =(X_{i},A_{i},P_{i},B_{i})_{i\in I}$\textit{\ \ be an
abstract economy, where }$I$\textit{\ is a (possibly uncountable) set of
agents such that for each }$i\in I:$
\end{theorem}

1)\textit{\ }$X_{i}$\textit{\ is a non-empty compact set in a Hausdorff
locally convex topological vector space }$E_{i}$\textit{;}

2)\textit{\ }$B_{i}$\textit{\ is w-lower semicontinuous with nonempty convex
values;}

3)\textit{\ the set }$W_{i}:$\textit{\ }$=\left\{ x\in X\text{ / }\left(
A_{i}\cap P_{i}\right) (x)\neq \emptyset \right\} $ \textit{is open;}

4)\textit{\ }$(A_{i}\cap P_{i})$\textit{\ has the e-LSCS-property on }$W_{i}$%
\textit{.}

\textit{Then there exists an equilibrium point }$x^{\ast }\in X$ \textit{\
for }$\Gamma $\textit{,}$\ i.e.$\textit{, for each }$i\in I$\textit{, }$%
x_{i}^{\ast }\in \overline{B}_{i}(x^{\ast })$\textit{\ and }$(A_{i}\cap
P_{i})(x^{\ast })=\emptyset $ for each $i\in I$\textit{.\medskip }

\textit{Proof.} For each\textit{\ }$i\in I$, let \ss $_{i}$ denote the
family of all open absolutely convex neighborhoods of zero in $E_{i}.$ Let $%
V=(V_{i})_{i\in I}\in \tprod\limits_{i\in I}$\ss $_{i}.$ Since $(A_{i}\cap
P_{i})$ has the e-LSCS-property on $W_{i}$, it follows that there exists a
lower semicontinuous\textit{\ }correspondence $F_{i}^{V_{i}}:X\rightarrow
2^{X_{i}}$ with convex values such that $F_{i}^{V_{i}}(x)\subset (A_{i}\cap
P_{i})(x)+V_{i}$ and $x_{i}\notin $cl$F_{i}^{V_{i}}(x)$ for each $x\in W_{i}$%
.

Let's define the correspondence $T_{i}^{V_{i}}:X\rightarrow 2^{X_{i}}$, by

$T_{i}^{V_{i}}(x):=\left\{ 
\begin{array}{c}
F_{i}^{V_{i}}(x)\text{, \ \ \ \ \ \ \ \ \ \ \ \ \ \ \ \ \ \ \ \ \ \ \ \ if }%
x\in \text{cl}W_{i}\text{, } \\ 
(B_{i}(x)+V_{i})\cap X_{i}\text{, \ \ \ \ \ \ \ \ \ \ if }x\notin \text{cl}%
W_{i}\text{.}%
\end{array}%
\right. $

$B^{V_{i}}:X\rightarrow 2^{X_{i}},$ $B^{V_{i}}(x)=(B_{i}(x)+V_{i})\cap X_{i}$
is lower semicontinuous according to the assumption 2).

Let $U$ be an closed subset of $\ X_{i}$, then

$U^{^{\prime }}:=\{x\in X$ $\mid T_{i}^{V_{i}}(x)\subset U\}$

\ \ \ =$\{x\in $cl$W_{i}\mid T_{i}^{V_{i}}(x)\subset U\}\cup \{x\in
X\setminus $cl$W_{i}\mid T_{i}^{V_{i}}(x)\subset U\}$

\ \ \ =$\left\{ x\in \text{cl}W_{i}\mid F_{i}^{V_{i}}(x)\subset U\right\}
\cup \left\{ x\in X\mid (B_{i}(x)+V_{i})\cap X_{i}\subset U\right\} $

\ \ \ 

$U^{^{\prime }}$ is a closed set, because cl$W_{i}$ is closed, $\left\{ x\in 
\text{cl}W_{i}\text{ }\mid F_{i}^{V_{i}}(x)\subset U\right\} $ is closed
since $F_{i}^{V_{i}}(x)$ is lower semicontinuous map on cl$W_{i}$ and the set

$\left\{ x\in X\mid \text{ }(B_{i}(x)+V_{i})\cap X_{i}\subset U\right\} $ is
closed since $(B_{i}(x)+V_{i})\cap X_{i}$ is lower semicontinuous. Then $%
T_{i}^{V_{i}}$ is lower semicontinuous on $X.$ By Theorem 7.3.3 in [8], cl$%
T_{i}^{V_{i}}$ is lower semicontinuous and has closed convex values.

Since $X$ is a compact convex set, by Wu's fixed-point theorem [18], there
exists $x_{V}^{\ast }\in X$ such that for each $i\in I$, $(x_{V}^{\ast
})_{i}\in $cl$T_{i}^{V_{i}}(x_{V}^{\ast })$. If $x_{V}^{\ast }\in $cl$W_{i},$
$(x_{V}^{\ast })_{i}\in $cl$F_{i}^{V_{i}}(x_{V}^{\ast })$, which is a
contradiction.

Hence, $(x_{V}^{\ast })_{i}\in $cl[$(B_{i}(x_{V}^{\ast })+V_{i})\cap X_{i}]$
and $(A_{i}\cap P_{i})(x_{V}^{\ast })=\emptyset .$ We have that $%
(x_{V}^{\ast })_{i}\in $cl[$(B_{i}(x_{V}^{\ast })+V_{i})\cap X_{i}]\subset 
\overline{B^{V_{i}}}(x_{V}^{\ast }).$

Let's define $Q_{V}=\cap _{i\in I}\{x\in X:$ $x_{i}\in \overline{B^{V_{i}}}%
(x)$ and $(A_{i}\cap P_{i})(x)=\emptyset \}.$ We have that $x_{V}^{\ast }\in
Q_{V},$ then $Q_{V}\neq \emptyset $ and so it is a non-empty closed subset
of $X$ by 3) and hence $Q_{V}$ is compact.

We prove that the family $\{Q_{V}:V\in \underset{i\in I}{\tprod }\text{\ss }%
_{i}\}$ has the finite intersection property.

Let $\{V^{(1)},V^{(2)},...,V^{(n)}\}$ be any finite set of $\underset{i\in I}%
{\tprod \text{\ss }_{i}}$ and let $V^{(k)}=\underset{i\in I}{\tprod }%
V_{i}^{(k)}$, where $V_{i}^{(k)}\in $\ss $_{i}$ for each $i\in I.$ Let $%
V_{i}=\underset{k=1}{\overset{n}{\cap }}V_{i}^{(k)}$, then $V_{i}\in \text{%
\ss }_{i};$ thus $V=\underset{i\in I}{\tprod }V_{i}\in \underset{i\in I}{%
\tprod }\text{\ss }_{i}.$ Clearly, $Q_{V}\subset \underset{k=1}{\overset{n}{%
\cap }}Q_{V^{(k)}}$ so that $\underset{k=1}{\overset{n}{\cap }}%
Q_{V^{(k)}}\neq \emptyset .$

Therefore, the family $\{Q_{V}:V\in \underset{i\in I}{\tprod }\text{\ss }%
_{i}\}$ has the finite intersection property. Since $X$ is compact, we have
that $\cap \{Q_{V}:V\in \underset{i\in I}{\tprod }\text{\ss }_{i}\}\neq
\emptyset .$ Let's take any $x^{\ast }\in \cap \{Q_{V}:V\in \underset{i\in I}%
{\tprod }\text{\ss }_{i}\},$ then for each $i\in I$ and each $V_{i}\in \text{%
\ss }_{i},$ $x_{i}^{\ast }\in \overline{B^{V_{i}}}(x^{\ast })$ and $%
(A_{i}\cap P_{i})(x^{\ast })=\emptyset ;$ but then, $x_{i}^{\ast }\in 
\overline{B}_{i}(x^{\ast })$ from Lemma 1 and $(A_{i}\cap P_{i})(x^{\ast
})=\emptyset $ for each $i\in I.$\medskip

\section{The model of a generalized multiobjective game and the existence of
generalized Pareto equilibrium}

The purpose of this section is to make a preliminary unitary presentation of
the model of a constrained multicriteria game in its strategic form and of
the solution concepts for this type of game, and also to state an existence
result for generalized Pareto equilibria.\medskip

Let $I$ be a finite set (the set of players). For each $i\in I,$ let $X_{i}$
be the set of strategies and define $X=\tprod\nolimits_{i\in I}X_{i}$. Let $%
T^{i}:X\rightarrow 2^{\mathbb{R}^{k_{i}}}$, where $k_{i}\in \mathbb{N}$, be
the multicriteria payoff function and let $A^{i}:X\rightarrow 2^{X_{i}}$ be
a constraint correspondence.

\begin{definition}[9]
The family $G=(X_{i},A^{i},T^{i})_{i\in I}$ is called a \textit{generalized
multicriteria (multiobjective) game}.\medskip
\end{definition}

Any n-tuple of strategies can be regarded as a point in the product space of
sets of players' strategies: $x=(x_{1},x_{2},...x_{n})\in X.$ For each
player $i\in I,$ the vector of the $n-1$ strategies of the other ones will
be denoted by $x_{-i}=(x_{1},...x_{i-1},x_{i+1},...x_{n})\in
X_{-i}=\tprod\nolimits_{j\in I\setminus \{i\}}X_{i}.$ We note that $%
x=(x_{-i},x_{i}).$

We assume that each player is trying to minimize his/her own payoff
according to his/her preferences, where for each player $i\in I,$ the
preference $"\gtrsim _{i}"$ over the outcome space $\mathbb{R}^{k_{i}}$ is
the following:

$z^{1}\gtrsim _{i}z^{2}$ if only if $z_{j}^{1}\geq z_{j}^{2}$ for each $%
j=1,2,...k_{i}$ and $z^{1},z^{2}\in \mathbb{R}^{k_{i}}.$ The following
preference can be defined on $X$ for each player $i$ (see [9]):

$x\gtrsim _{i}y$ whenever $F^{i}(x)\gtrsim _{i}F^{i}(y)$ and $x,y\in X.$

Let $x^{\ast }=(x_{1}^{\ast },x_{2}^{\ast },...,x_{n}^{\ast })\in X.$

We introduce slight generalizations of the equilibrium concepts defined by
Kim and Ding in [9].

\begin{definition}
A strategy $x_{i}^{\ast }\in X_{i}$ of player $i$ is said to be a \textit{%
generalized} \textit{Pareto} \textit{efficient strategy (respectively, a
weak Pareto efficient strategy) with respect to }$x$ if $x_{i}^{\ast }\in 
\overline{A^{i}}(x^{\ast })$ and there is no strategy $x_{i}\in
A^{i}(x^{\ast })$ such that
\end{definition}

$T^{i}(x^{\ast })-T^{i}(x_{i-1}^{\ast },x_{i})\in \mathbb{R}%
_{+}^{k_{i}}\backslash \{0\}$ (respectively, $T^{i}(x^{\ast
})-T^{i}(x_{i-1}^{\ast },x_{i})\in $int$\mathbb{R}_{+}^{k_{i}}\backslash
\{0\}).$

\begin{definition}
A strategy $x^{\ast }\in X$ is said to be a \textit{generalized Pareto} 
\textit{equilibrium (respectively, a weak Pareto equilibrium) of a game }$%
G=(X_{i},A^{i},T^{i})_{i\in I},$\textit{\ }if for each player $i\in I$, $%
x_{i}^{\ast }\in X_{i}$ is a Pareto efficient strategy against $x^{\ast }$
(respective, a generalized weak Pareto efficient strategy against $x^{\ast }$%
)$.$
\end{definition}

The following notion contains the idea of a game equilibrium defined by
using a scalarization function. In this case, the scalarization method uses
weighted coefficients $W_{i},$ so that each player $i$ has his own vector of
weights $W_{i}\in \mathbb{R}_{+}^{k_{i}}\backslash \{0\}.$

\begin{definition}
A strategy $x^{\ast }\in X$ is said to be a \textit{generalized weighted Nash%
} \textit{equilibrium} with respect to the weighted vector $W=(W_{i})_{i\in
I}$ with $W_{i}=(W_{i,1},W_{i,2},$
\end{definition}

\noindent $...,W_{i,k_{i}})\in \mathbb{R}_{+}^{k_{i}}$ of the multiobjective
game $G=((X_{i},A^{i},T^{i})_{i\in I}$ if for each player $i\in I$, we have

1) $x_{i}^{\ast }\in \overline{A^{i}}(x^{\ast });$

2) $W_{i}\in \mathbb{R}_{+}^{k_{i}}\backslash \{0\};$

3) for all $x_{i}\in A^{i}(x^{\ast }),$ $W_{i}\cdot T^{i}(x^{\ast })\leq
W_{i}\cdot T^{i}(x_{-i}^{\ast },x_{i}),$ where $\cdot $ denotes the inner
product in $\mathbb{R}^{k_{i}}.$

\begin{remark}
In particular, if $W_{i}\in \Delta ^{k_{i}}=\{u_{i}\in \mathbb{R}%
_{+}^{k_{i}} $ with $\tsum\nolimits_{j=1}^{k_{i}}u_{i,j}=1\}$ for each $i\in
I,$ then the strategy $x^{\ast }\in X$ is said to be a normalized
generalized weighted Nash equilibrium with respect to $W.$
\end{remark}

\begin{remark}
If for each $i\in I,$ $\overline{A^{i}}$ has closed values and a closed
graph in $X\times X_{i},$ the notions of equilibrium introduced above
coincide with the equilibrium notions defined by Kim and Ding in [9].
\end{remark}

The relationship between the two types of equilibrium notions is given by
the following result.

\begin{lemma}
Each normalized generalized weighted Nash equilibrium $x^{\ast }\in X$ with
a weight $W=(W_{1},...W_{n})\in \Delta ^{k_{1}}\times ...\times \Delta
^{k_{n}}$ (respectively, $W=(W_{1},...,W_{n})\in $ int$\Delta ^{k_{1}}\times
...\times $int$\Delta ^{k_{n}}$) is a weak Pareto equilibrium (respectively,
a Pareto equilibrium) of the game $G=((X_{i},A^{i},T^{i})_{i\in I}.$
\end{lemma}

The proof follows the same line as the proof of Lemma 7 in [9].

\begin{remark}
As in [9], the above lemma remains true when $W=(W_{1},...,W_{n})$ satisfies 
$W_{i}\in \mathbb{R}_{+}^{k_{i}}$ (resp. $W_{i}\in $int$\mathbb{R}%
_{+}^{k_{i}}).$
\end{remark}

In order to prove the existence result for generalized weighted Nash
equilibrium of generalized multiobjective games, first we prove the
following lemma.

\begin{lemma}
Let $X$ be a nonempty convex compact of a Hausdorff locally convex
topological vector space $E,$ $D$ be a nonempty compact convex and
metrizable subset of $X$, $A:X\rightarrow 2^{X}$ be a correspondence with
non-empty convex values and $f:X\times X\rightarrow \overline{\mathbb{R}}$
be a function such that:
\end{lemma}

1)\textit{\ }$A$\textit{\ is almost weakly lower semicontinuous with respect
to }$D$\textit{\ and }$\overline{A^{V}}$\textit{\ is convex nonempty valued
for each open absolutely convex symmetric neighbourhood }$V$\textit{\ of }$0$%
\textit{\ in }$E$\textit{;}

2)\textit{\ The correspondence }$F:X\rightarrow 2^{X}$\textit{, }$%
F(x)=\{y\in X:f(x,x)-f(y,x)>0\}$\textit{\ is almost weakly lower
semicontinuous with respect to }$D$\textit{\ on }$K=\{x\in X:x\in \overline{A%
}(x)\}$ \textit{and }$\overline{F^{V}}$\textit{\ is convex valued for each
open absolutely convex symmetric neighbourhood }$V$\textit{\ of }$0$\textit{%
\ in }$E$\textit{;}

3)\textit{\ }$x\notin \overline{F}(x)$\textit{\ for each }$x\in K$\textit{;}

\textit{then there exists }$x^{\ast }\in X$\textit{\ such that }$x^{\ast
}\in \overline{A}(x^{\ast })$\textit{\ and }$f(x^{\ast },x^{\ast })\leq
f(y,x^{\ast })$\textit{\ for each }$y\in \overline{A}(x^{\ast }).$

\textit{Proof.} We notice first that\textit{\ }the set $K=\{x\in X:x\in 
\overline{A}(x)\}$\ is closed.

Assume that for each $x\in K,$ $A(x)\cap F(x)\neq \phi $ and define the
correspondence $G:X\rightarrow 2^{X}$ by

$G(x)=\left\{ 
\begin{array}{c}
A(x)\cap F(x)\text{ if }x\in K; \\ 
A(x)\text{ \ \ \ \ \ if \ \ \ \ }x\notin K.%
\end{array}%
\right. $

By 1) and 3), the correspondence $\overline{G^{V}}:X\rightarrow 2^{X}$ is
lower semicontinuous for each open absolutely convex symmetric neighbourhood 
$V$ of $0$ in $E,$ and has nonempty convex closed values (we can prove this
fact by using an argument similar with that one from the Theorem 3). By
Corollary 2, there exists $x^{\ast }\in X$ such that $x^{\ast }\in \overline{%
G}(x^{\ast }).$ By definition of $G$ and $A,$ $x^{\ast }$ must be in $K.$ It
follows that $x^{\ast }\in \overline{A\cap F}(x^{\ast })$, and since clGr$%
A\cap F\subset $clGr$A\cap $clGr$F$, we have that $x^{\ast }\in \overline{A}%
(x^{\ast })\cap \overline{F}(x^{\ast })$, that is $x^{\ast }\in \overline{F}%
(x^{\ast }),$ which contradicts 3). Therefore, there exists $x^{\ast }\in K$
such that $\overline{A}(x^{\ast })\cap \overline{F}(x^{\ast })=\phi $ (this
implies also $A(x^{\ast })\cap F(x^{\ast })=\phi )$. Hence

\begin{center}
$x^{\ast }\in \overline{A}(x^{\ast })$ and $f(x^{\ast },x^{\ast })\leq
f(y,x^{\ast })$ for each $y\in \overline{A}(x^{\ast }).$
\end{center}

\begin{example}
Let $f:[-1,1]\times \lbrack -1,1]\rightarrow \mathbb{R},$
\end{example}

$f(x,y)=\left\{ 
\begin{array}{c}
1\text{ \ \ \ \ \ \ \ \ \ \ \ \ \ \ \ \ \ \ \ \ \ \ \ \ \ \ \ \ \ \ \ if \ \
\ \ \ \ \ \ \ \ }(x,y)=(-1,0); \\ 
2\text{ \ \ \ \ \ \ \ \ \ \ \ \ \ \ \ \ \ \ \ \ \ \ \ \ \ \ \ \ \ \ \ \ \ if
\ \ \ \ \ \ \ \ \ \ }(x,y)=(0,0); \\ 
1\text{ \ \ \ \ \ \ \ \ \ \ \ \ \ \ \ \ \ \ \ \ \ \ \ \ if \ }x,y\in \lbrack
0,1)\times \lbrack 0,1]\backslash \{(0,0)\}; \\ 
2\text{ \ \ \ \ \ \ \ \ \ \ \ \ \ \ \ \ \ \ if \ \ \ \ \ \ \ \ \ \ \ \ \ }%
(x,y)\in (\frac{1}{2},1]\times \lbrack -1,0); \\ 
3\text{\ \ \ \ if \ \ \ }(x,y)\in \lbrack -1,\frac{1}{2}]\times \lbrack
-1,0)\cup \{(-1,0)\times \{0\}\}; \\ 
4\text{ \ \ \ \ \ \ \ \ \ \ \ \ \ \ \ \ \ if \ \ \ \ \ \ \ \ \ \ \ \ \ \ \ }%
(x,y)\in \lbrack -1,0)\times (0,1]; \\ 
0\text{ \ \ \ \ \ \ \ \ \ \ \ \ \ \ \ \ \ \ \ \ if \ \ \ \ \ \ \ \ \ \ \ \ \
\ \ \ }(x,y)\in \{1\}\times (0,1];%
\end{array}%
\right. $

Let $A:[-1,1]\rightarrow 2^{[-1,1]}$ defined by $A(x)=[-1,0]$ if $x\in
\lbrack -1,1].$

$A$ is lower semicontiuous on $[-1,1]$ and $K=\{x\in \lbrack -1,1]:x\in 
\overline{A}(x)\}=[-1,0]$ is closed.

$F:X\rightarrow 2^{X}$, $F(x)=\{y\in X:f(x,x)-f(y,x)>0\}$

$\ \ \ \ \ \ \ \ \ \ \ \ \ \ \ \ \ \ \ \ \ \ \ \ \ \ F(x)=\left\{ 
\begin{array}{c}
(\frac{1}{2},1]\text{ if }x\in \lbrack -1,0); \\ 
\{-1\}\text{ \ \ \ if \ \ \ \ }x=0; \\ 
\{1\}\text{ \ \ if \ \ }x\in (0,1].%
\end{array}%
\right. $

$F$ is not lower semicontinuous and $x\notin \overline{F}(x),$ $\forall x\in
K=[-1,0]$, where

$\overline{F}_{K}(x)=\left\{ 
\begin{array}{c}
\lbrack \frac{1}{2},1]\text{ \ \ \ \ \ if \ \ \ \ }x\in \lbrack -1,0); \\ 
\{-1\}\cup \lbrack \frac{1}{2},1]\text{ \ \ \ if \ \ \ \ }x=0;%
\end{array}%
\right. $

We also have that for each $V=(-\varepsilon ,\varepsilon )$ with $%
\varepsilon >0,$ and $D=[0,1],$ $F_{|K}^{V}$ is lower semicontinuous.

Therefore,

$F_{|K}(x)+(-\varepsilon ,\varepsilon )=\left\{ 
\begin{array}{c}
(\frac{1}{2}-\varepsilon ,1+\varepsilon )\text{ if }x\in \lbrack -1,0); \\ 
(-1-\varepsilon ,-1+\varepsilon )\text{\ \ if \ \ }x=0;%
\end{array}%
\right. $

For $\varepsilon \in (0,\frac{1}{2}],$ $F_{|K}^{V}(x)=\left\{ 
\begin{array}{c}
(\frac{1}{2}-\varepsilon ,1]\text{ if }x\in \lbrack -1,0); \\ 
\phi \text{ \ \ \ \ \ \ \ if \ \ \ \ \ \ \ \ }x=0;%
\end{array}%
\right. $

For $\varepsilon \in (\frac{1}{2},1],$ $F_{|K}^{V}(x)=\left\{ 
\begin{array}{c}
\lbrack 0,1]\text{ \ \ if }x\in \lbrack -1,0); \\ 
\phi \text{\ \ \ \ \ \ \ \ \ if \ \ \ \ \ \ }x=0;%
\end{array}%
\right. $

For $\varepsilon \in (1,2],$ $F_{|K}^{V}(x)=\left\{ 
\begin{array}{c}
\lbrack 0,1]\text{ if \ }x\in \lbrack -1,0); \\ 
\lbrack 0,-1+\varepsilon )\text{\ if\ }x=0;%
\end{array}%
\right. $

For $\varepsilon >2,$ $F_{|K}^{V}(x)=[0,1]$, $x\in \lbrack -1,0].$

Then,

For $\varepsilon \in (0,\frac{1}{2}],$ $\overline{F_{|K}^{V}}(x)=[\frac{1}{2}%
-\varepsilon ,1]$ if $x\in \lbrack -1,0];$

For $\varepsilon >\frac{1}{2}$ $\overline{F_{|K}^{V}}(x)=[0,1]$ \ \ if \ $%
x\in \lbrack -1,0];$

$\overline{F_{|K}^{V}}$ is lower semicontinuous and nonempty convex valued.

By Lemma 3, we have that there is $x^{\ast }\in \overline{A}(x^{\ast })$
such that $A(x^{\ast })\cap F(x^{\ast })=\phi .$

For example, $x^{\ast }=-\frac{1}{2},$ -$\frac{1}{2}\in \overline{A}(-\frac{1%
}{2})$ and -$\frac{1}{2}\notin F(-\frac{1}{2}),$ that is $3=f(-\frac{1}{2},-%
\frac{1}{2})\geq f(y,-\frac{1}{2})=3$ for each $y\in \overline{A}(-\frac{1}{2%
})=[-1,0].$

Now, as an application of Lemma 3, we have the following existence theorem
of generalized weighted Nash equilibrium for generalized multiobjective
games.

\begin{theorem}
Let $I$ be a finite set of indices, let ($X_{i},A^{i},T^{i}$)$_{i\in I}$ be
a constrained multi-criteria game with for each $i\in I,$ $X_{i}$ is a
nonempty convex subset of a Hausdorff locally convex topological vector
space $E^{i}$ and suppose that there is a nonempty compact convex and
metrizable subset $D$ of $X=\tprod\nolimits_{i\in I}X_{i}$ and a weighted
vector $W=(W_{1},W_{2},...,W_{n})$ with $W_{i}\in \mathbb{R}%
_{+}^{k_{i}}\backslash \{0\}$ such that the following conditions are
satisfied:
\end{theorem}

1)\textit{\ for each }$i\in I,$\textit{\ }$A^{i}$\textit{\ is almost weakly
lower semicontinuous with respect to }$D$\textit{\ and }$\overline{%
A^{i,V_{i}}}$\textit{\ is convex nonempty valued for each open absolutely
convex symmetric neighbourhood }$V_{i}$\textit{\ of }$0$\textit{\ in }$E_{i}$%
\textit{;}

2)\textit{\ The set }$K=\{x\in X:x\in \overline{A}(x)\}$\textit{, where }$%
A(x)=\tprod\nolimits_{i\in I}A^{i}(x)$\textit{, is closed in }$X$\textit{;}

3)\textit{\ The correspondence }$F:X\rightarrow 2^{X}$\textit{, }$%
F(x)=\{y\in X:\tsum\nolimits_{i=1}^{n}W_{i}\cdot
(T^{i}(x_{-i},x_{i})-T^{i}(x_{-i},y_{i}))>0\}$\textit{\ is almost weakly
lower semicontinuous with respect to }$D$\textit{\ on }$K$ \textit{and }$%
\overline{F^{V}}$\textit{\ is convex valued for each open absolutely convex
symmetric neighbourhood }$V$\textit{\ of }$0$\textit{\ in }$E$\textit{;}

4)\textit{\ }$x\notin \overline{F}(x)$\textit{\ for each }$x\in K$\textit{;}

\textit{then there exists }$x^{\ast }\in X$\textit{\ such that }$x^{\ast }$%
\textit{\ is a generalized weighted Nash equilibria with respect to }$%
W.\medskip $

\textit{Proof.} Define the function $f:X\times X\rightarrow \mathbb{R}$ by $%
f(x,y)=\tsum\nolimits_{i=1}^{n}W_{i}\cdot
(T^{i}(x_{-i},x_{i})-T^{i}(x_{-i},y_{i})),$ $(x,y)\in X\times X.$ $\ $\ It
is easy to see that $f$ satisfies all hypothesis of Lemma 3, hence there
exists $x^{\ast }\in X$ such that $x^{\ast }\in \overline{A}(x^{\ast })$ and 
$\tsum\nolimits_{i=1}^{n}W_{i}\cdot (T^{i}(x_{-i}^{\ast },x_{i}^{\ast
})-T^{i}(x_{-i}^{\ast },y_{i})\leq 0$ for any $y\in \overline{A}(x^{\ast }).$
We use the fact that $\tprod\nolimits_{i\in I}A^{i}\subseteq \overline{%
\tprod\nolimits_{i\in I}A^{i}}\subseteq \tprod\nolimits_{i\in I}\overline{%
A^{i}}.$ We obtain first $x_{i}^{\ast }\in \overline{A_{i}}(x^{\ast })$ for
each $i\in I$. For any given $i\in I$ and any given $y_{i}\in A^{i}(x^{\ast
}),$ let $y=(x_{-i}^{\ast },y_{i}).$ Then

$W_{i}\cdot (T^{i}(x_{-i}^{\ast },x_{i}^{\ast })-T^{i}(x_{-i}^{\ast
},y_{i}))=$

$=\tsum\nolimits_{j=1}^{n}W_{j}\cdot (T^{j}(x_{-i}^{\ast },x_{i}^{\ast
})-T^{i}(x_{-i}^{\ast },y_{i}))-\tsum\nolimits_{j\neq i}W_{j}\cdot
(T^{j}(x_{-i}^{\ast },x_{i}^{\ast })-T^{i}(x_{-i}^{\ast },y_{i}))$

$=\tsum\nolimits_{j=1}^{n}W_{j}\cdot (T^{j}(x_{-i}^{\ast },x_{i}^{\ast
})-T^{i}(x_{-i}^{\ast },y_{i}))\leq 0.$

Therefore, we have $W_{i}\cdot (T^{i}(x_{-i}^{\ast },x_{i}^{\ast
})-T^{i}(x_{-i}^{\ast },y_{i}))$ $\leq 0$ for each $i\in I$ and $y_{i}\in
A^{i}(x^{\ast }).$ Hence, $x^{\ast }$ is a generalized weighted Nash
equilibrium of the game $G$ with respect to $W.\medskip $

By using Lemma 3, we obtain the following existence theorem of generalized
Pareto equilibrium as a consequence of Theorem 6.

\begin{theorem}
Let $I$ be a finite set of indices, let ($X_{i},A^{i},T^{i}$)$_{i\in I}$ be
a constrained multi-criteria game with for each $i\in I,$ $X_{i}$ is a
nonempty convex subset of a Hausdorff locally convex topological vector
space $E^{i}$ and suppose that there is a nonempty compact convex and
metrizable subset $D$ of $X=\tprod\nolimits_{i\in I}X_{i}$ and a weighted
vector $W=(W_{1},W_{2},...,W_{n})$ with $W_{i}\in \mathbb{R}%
_{+}^{k_{i}}\backslash \{0\}$ such that the following conditions are
satisfied:
\end{theorem}

1)\textit{\ for each }$i\in I,$\textit{\ }$A^{i}$\textit{\ is almost weakly
lower semicontinuous with respect to }$D$\textit{\ and }$\overline{%
A^{i,V_{i}}}$\textit{\ is convex nonempty valued for each open absolutely
convex symmetric neighbourhood }$V_{i}$\textit{\ of }$0$\textit{\ in }$E_{i}$%
\textit{;}

2)\textit{\ The set }$K=\{x\in X:x\in \overline{A}(x)\}$\textit{, where }$%
A(x)=\tprod\nolimits_{i\in I}A^{i}(x)$\textit{, is closed in }$X$\textit{;}

3)\textit{\ The correspondence }$F:X\rightarrow 2^{X}$\textit{, }$%
F(x)=\{y\in X:\tsum\nolimits_{i=1}^{n}W_{i}\cdot
(T^{i}(x_{-i},x_{i})-T^{i}(x_{-i},y_{i}))>0\}$\textit{\ is almost weakly
lower semicontinuous with respect to }$D$\textit{\ on }$K$ \textit{and }$%
\overline{F^{V}}$\textit{\ is convex valued for each open absolutely convex
symmetric }$V$\textit{\ of }$0$\textit{\ in }$E$\textit{;}

4)\textit{\ }$x\notin \overline{F}(x)$\textit{\ for each }$x\in K$\textit{;}

\textit{then there exists }$x^{\ast }\in X$\textit{\ such that }$x^{\ast }$%
\textit{\ is a generalized weak Pareto equilibrium.}

\textit{Furthermore, if }$W_{i}\in $\textit{int}$R_{+}^{k_{i}}\backslash
\{0\}$\textit{\ for all }$i\in I,$\textit{\ then }$x^{\ast }$\textit{\ is a
generalized Pareto equilibrium.}$\medskip $

The author thanks to Professor Jo\~{a}o Paulo Costa from the University of
Coimbra for the fruitfull discussions and for the hospitality he proved
during the visit to his departament.

\end{document}